\documentclass[11pt]{article}
\usepackage{amsmath,amsfonts,amssymb,mathrsfs,cases,mathdots,colordvi,url,extarrows,graphicx}

\usepackage{relsize}
\usepackage{mathrsfs}
\usepackage{color}
\usepackage{indentfirst}
\usepackage{float}

\newtheorem{defn}{Definition}[section]

\newtheorem{lemma}{Lemma}[section]
\newtheorem{thm}{Theorem}[section]
\newtheorem{prop}{Proposition}[section]

\newtheorem{example}{Example}[section]

\renewcommand{\Box}{\rule{2.2mm}{2.2mm}}

\def\beginproof{\par\noindent {\bf Proof.}\ \ }
\def\endproof{\hskip .5cm $\Box$ \vskip .5cm}

\def\beginproof{\par\noindent {\bf Proof.}\ \ }
\def\endproof{\hskip .5cm $\Box$ \vskip .5cm}

\begin{document}
\title{Directional necessary optimality conditions for bilevel programs}
\author{Kuang Bai\thanks{Department of Applied Mathematics, The Hong Kong Polytechnic University, Hong Kong, China.
Email: kuang.bai@polyu.edu.hk.} \and Jane J. Ye \thanks{Corresponding author. Department of Mathematics and Statistics, University of Victoria, Canada. The research of this author was partially
supported by NSERC. Email: janeye@uvic.ca.}}
  \date{}
\maketitle 
\begin{abstract} 
	The bilevel program  is an optimization problem where the constraint  involves  solutions to a parametric optimization problem. It is well-known that the value function reformulation provides an equivalent single-level optimization problem  but it results in   a nonsmooth optimization problem which never satisfies the usual constraint qualification such as the  Mangasarian-Fromovitz constraint qualification (MFCQ). In this paper we show that even the first order sufficient condition for metric subregularity  (which is in general weaker than  MFCQ) fails at each feasible point of the bilevel program.  We introduce the concept of directional calmness condition and show that under {the} directional calmness condition, the directional necessary optimality condition  holds. {While the directional optimality condition is in general sharper than the non-directional one,}  the directional calmness condition is in general weaker than the classical calmness condition and hence is more likely to hold. {We perform the directional sensitivity analysis of the value function and} propose the directional quasi-normality as a sufficient condition for the directional calmness. An example is given to show that the directional quasi-normality condition may hold for the bilevel program.

	\vskip 10 true pt
	
	\noindent {\bf Key words.}\quad bilevel programs,  constraint qualifications, necessary optimality conditions,  directional derivatives,  directional subdifferentials, directional quasi-normality

	\vskip 10 true pt
	
	\noindent {\bf AMS subject classification:} {90C30, 91A65,  49K40}.
	
\end{abstract}

\section{Introduction}

The motivation for studying bilevel optimization originated in economics under the name of Stackelberg games \cite{Stackelberg} since 1934. In economics, it is used to model  interactions between a leader and its follower of a two level hierarchical system and hence is referred to as leader and follower games or principal-agent problems. In recent years, bilevel programs find wider range of applications (see e.g. \cite{Dempe, Bireview, Luo, Outrata} and   references within).
In particular, bilevel programs have been used to model hyper-parameter selection in machine learning (see e.g. \cite{ML,tax}) in recent years.

In this paper, we consider bilevel programs in the following form:
\begin{align*}
\mbox{(BP)}\quad &\min_{x,y} \quad F(x,y)\\
&s.t.\quad y\in S(x),\ G(x,y)\leq0,
\end{align*}
where for any given $x$, $S(x)$ denotes the solution set of the lower level program
\[
(P_x)\quad \min_{y} f(x,y) \quad \mbox{s.t. }  g(x,y)\leq 0,
\] 
and $F, f:\mathbb R^n\times\mathbb R^m\rightarrow\mathbb R,\ G:\mathbb R^n\times\mathbb R^m\rightarrow\mathbb R^q$, $ g:\mathbb R^n\times\mathbb R^m\rightarrow\mathbb R^p$ are continuously differentiable.


To obtain an optimality condition for (BP), one may reformulate it as a single-level optimization problem and apply optimality conditions to the single-level problem. There are three approaches for reformulating (BP)  as a single-level optimization problem in the literature. The earliest approach is the so-called first order approach or the Karush-Kuhn-Tucker (KKT) approach by which one  replaces the constraint $y\in S(x)$ by its first order optimality conditions and minimizing over the original variables as well as the multipliers. The resulting single-level optimization problem is the
so-called mathematical program with equilibrium constraints (MPEC),  which was popularly studied over the last three decades; see e.g. \cite{Luo,Outrata} for the general theory  and  \cite{Ye00,YY,GY19} for the optimality conditions derived by using this approach. 
The value function approach first proposed in \cite{VP}
  replaces the constraint $y\in S(x)$ by $f(x,y)-V(x)\leq 0, g(x,y)\leq 0$, where $V(x):=\inf_{y}\{f(x,y)| g(x,y)\leq0\}$ is the value function of the lower level program $(P_x)$. And the combined approach (\cite{YZ, Ye11}) not only replaces the constraint $y\in S(x)$ by $f(x,y)-V(x)\leq 0, g(x,y)\leq 0$ but also adds the first order optimality conditions. The first order approach is obviously only applicable if the first order optimality condition is necessary and sufficient for optimality; e.g. when the lower level program is convex and certain constraint qualification holds. Both the KKT approach and the combined approach suffer from the drawback that the resulting MPEC may not be equivalent to the original (BP) if the local optimality is considered; see \cite{DD} for the discussion for the KKT approach and \cite{Yebook} for the combined approach.

In this paper by the value function approach, we reformulate (BP) as the following equivalent problem:
\begin{align*}
\mbox{(VP)}~~~~\quad \min_{x,y}\quad &F(x,y)\\
{\rm s.t.}\quad &f(x,y)-V(x)\leq0, \ g(x,y)\leq0, \
G(x,y)\leq0.
\end{align*}
Under fairly reasonable assumptions, the value function $V(x)$ is Lipschitz continuous and hence a  nonsmooth Fritz John type necessary optimality condition holds at a local optimal solution. For {a} KKT type necessary optimality condition to hold, in general one needs to assume certain constraint qualifications. Unfortunately, it is known 
 (\cite[Proposition 3.1]{YZ95}) that   the nonsmooth MFCQ or equivalently the no nonzero abnormal multiplier constraint qualification (NNAMCQ), a standard constraint qualification for  nonsmooth mathematical programs,  fails to hold at any feasible point of (VP). 
For an optimization problem with Lipschitz continuous problem data, it is known that the necessary optimality condition  holds provided that the problem is calm  in the sense of Clarke \cite[Definition 6.41]{Clarke}. Ye and {Zhu} \cite{YZ95} introduced the partial calmness condition for problem (VP) which means that a local solution of problem (VP) is also a local solution of the partially penalized problem for certain $\rho>0$
\begin{eqnarray*}
\mbox{(VP)}_\rho~~~~	\min_{x,y} && F(x,y)+\rho (f(x,y)-V(x))\\
	s.t. && g(x,y)\leq 0, G(x,y)\leq 0.
	\end{eqnarray*}  
Since the most difficult constraint $f(x,y)-V(x)$ is moved to the objective,  the KKT condition would hold under some  constraint qualifications for the partially penalized problem $\mbox{(VP)}_\rho$. It is easy to show that the {full} calmness implies the partial calmness and the partial calmness plus the {full} calmness of the partially penalized problem $\mbox{(VP)}_\rho$ implies the {full} calmness condition for problem (VP).
Some sufficient conditions for partial calmness and its relationship with exact penalization were further discussed in \cite{YZ95,YZ97,YZZ}.  Unfortunately for problem (VP), the partial calmness or the {full} calmness condition is still a fairly strong condition (see e.g., \cite{HS} for discussions).  {And there are very few constraint qualifications or sufficient conditions for partial calmness for (VP) in the literature. However the  single-level program  obtained from the combined approach is much more likely to satisfy the partial calmness condition. In fact, a recent paper \cite{KYYZ} showed that the partial calmness condition holds generically for the combined program when the upper level variable is one-dimensional. Recently, \cite{XY} has extended  the relaxed constant positive linear dependence constraint qualification (RCPLD) to  bilevel programs and has shown that it is a constraint qualification. } 

Recently  Gfrerer \cite[Theorem 7]{Gfr13} derived a directional version of the KKT type necessary optimality condition for  mathematical programs with a generalized equation constraint induced by a set-valued map under the directional metric subregularity constraint qualification. 
The directional KKT condition is in general sharper than the nondirectional KKT condition and the directional metric subregularity  is weaker than the nondirectional one. Inspired by this approach, in this paper we aim at developing a directional KKT condition for problem (VP). 
First we review the following  concept of directional neighborhood recently introduced by Gfrerer in \cite{Gfr13}.
  Given a direction $d\in \mathbb{R}^n$, and positive numbers $\epsilon,\delta>0$, the directional neighborhood of direction $d$ is a set defined by 
\begin{equation*}
{\cal V}_{\epsilon,\delta}(d):=\{z\in\epsilon\mathbb{B}|\big\| \| d\| z-\| z\| d\big\|\leq\delta \|z\| \|d\|\}.
\end{equation*}
It is easy to see that the directional neighborhood of direction $d=0$ is just the open ball $\epsilon\mathbb{B}$ and the directional neighborhood of a nonzero direction $d\not =0$ is a smaller subset of $\epsilon\mathbb{B}$.
Hence many regularity conditions can be extended to a directional version which is  weaker than the original nondirectional one. 
		We say that (VP) is calm at a feasible solution $(\bar x,\bar y)$ in direction $d\in \mathbb{R}^{n+m}$ if there exist positive scalars $\epsilon,\delta,\rho$, such that for any $\alpha\in \epsilon\mathbb B$ and any $(x,y)\in (\bar x,\bar y)+{\cal V}_{\epsilon,\delta}(d)$ satisfying {$\varphi(x,y)+\alpha\leq0$} with $\varphi(x,y):=(f(x,y)-V(x), g(x,y), G(x,y))$ one has,
		\begin{equation*}
	F(x,y)-	F(\bar x,\bar y)+\rho\|\alpha\|\geq0.
		\end{equation*}
	It is obvious that when the direction $d=0$, the directional calmness is reduced to the classical calmness condition \cite[Definition 6.41]{Clarke}. When $d\not =0$, since the directional neighborhood is in general smaller than the usual neighborhood, the directional calmness condition is in general weaker than the nondirectional calmness condition.
	It is obvious that if $(\bar x,\bar y)$ solves (VP), then under the calmness condition in direction $d$, $(\bar x,\bar y)$ is also a solution of the following penalized problem
	\begin{eqnarray*}
{\rm (DP)}~~~~	\min_{x,y} && F(x,y)+{\rho \|\varphi_+(x,y)\|}\\
	s.t. && (x,y)\in  (\bar x,\bar y)+{\cal V}_{\epsilon,\delta}(d).
	\end{eqnarray*}  
	The directionally penalized problem (DP) is much easier to deal  with than (VP) since all the inequality constraints are moved to the objective function. By using the nonsmooth calculus, one can then show that  $(\bar x,\bar y)$ satisfies a KKT condition provided the value function is Lipschitz continuous. In fact we can achieve more. {When $d$ is a critical direction,} we can show that $(\bar x,\bar y)$ satisfies a directional KKT condition in which a directional Clarke subdifferential (see Definition \ref{Clarke}) of the value function $V(x)$  at {$\bar x$} in direction $d$ is used instead of the Clarke subdifferential. Since the directional Clarke subdifferential is a subset of the Clarke subdifferential, the directional KKT condition is sharper than the nondirectional one.   The idea of deriving optimality conditions with respect to directions can also be found in \cite{DP}, where an optimality condition is formulated using the so-called directional convexificators. Note that similarly, the directional KKT condition can be extended to the combined program. But for the sake of simplicity and readability, we will leave this for future work.
	
	
	To make the directional calmness condition and the directional KKT condition useful, we have two issues to consider. 
	First, under what conditions, the value function is directionally Lipschitz continuous and directionally differentiable and how to calculate the directional limiting subdifferential and the directional derivative of the value function which will be needed in the directional KKT condition for problem (VP). In this paper, we have derived some formulas for the directional derivative of the value function and an upper  estimate for the Clarke directional subdifferential of the value function $V(x)$.
Secondly, how to derive a verifiable constraint qualification which ensures the directional calmness condition of (VP)? 
It is known that the first order sufficient condition for metric subregularity (FOSCMS) (introduced in Gfrerer and Klatte \cite[Corollary ]{GKlatte16} for the smooth case and \cite[Proposition 2.2]{BYZ} for the nonsmooth case)  is a sufficient condition for the metric subregularity of the set-valued map $\Phi(x,y):=\varphi(x,y)-\mathbb{R}^{p+q+1}_-$ which in turn implies the calmness of the problem (VP). FOSCMS is in general weaker than NNAMCQ and hence it is natural to ask if FOSCMS would hold for (VP).  Unfortunately in Proposition \ref{abnorm}, we show that FOSCMS also fails for problem (VP) in any critical direction. We propose the directional quasi-normality as a sufficient condition for the directional calmness condition and give an example to show that the directional quasi-normality is possible to hold for (VP).

	Other than deriving a weaker constraint qualification and a shaper necessary optimality condition for bilevel programs, we have also made contributions that are of independent interest as summarized below.
	\begin{itemize}
		\item We introduce the concept of directional Clarke subdifferentials and derive some useful calculus rules for directional subdifferentials; see Proposition \ref{Prop2.4}. 
		\item For an optimization problem with directionally Lipschitz continuous objective function and directionally Lipschitz  and directionally differentiable inequality  constraints, we derive a directional KKT condition under the  directional calmness condition; see Theorem \ref{dKKT}. An example of a bilevel program is given to show that the directional calmness is weaker than the classical calmness; see Example \ref{Ex3.1}. 
		 \item The classical results for the directional derivative of the value function are improved with weaker assumptions: see Propositions \ref{ddv} and \ref{Semic}. Sufficient conditions for directional Lipschitz continuity of the value function is given in Theorem \ref{Lip} and the upper estimate of the directional subdifferential of the value function is given in Theorems \ref{estimates} and \ref{subdiff}.
\end{itemize}

We organize the paper as follows. In the next section, we provide the notations, preliminaries and preliminary results. In Section 3 we derive the directional KKT condition under the directional calmness condition for a general optimization problem with directionally Lipschitz inequality constraints. In section 4, we study directional sensitivity analysis of the value function.  Finally in section 5, we apply the previous results to (VP) and derive a verifiable constraint qualification and a necessary optimality condition.

\section{Preliminaries} We first give notations that will be used in the paper. We denote by $\overline{\mathbb R}:=\mathbb R\cup\{\pm\infty\}$, while $\mathbb R:=(-\infty,+\infty)$. $\|\cdot\|$ denotes the Euclidean norm. $\langle a,b\rangle$ denotes the inner product of vectors $a,b$. Let $\Omega$ be a set. By $x^k\xrightarrow{\Omega}\bar{x}$ we mean $x^k\rightarrow\bar{x}$ and for each $k$, $x^k\in \Omega$.  By $x^k\xrightarrow{u}\bar x$ where $u$ is a vector, we mean that the sequence $\{x^k\}$ approaches $\bar x$ in direction $u$, i.e., there exist $t_k\downarrow 0, u^k\rightarrow u$ such that $x^k=\bar x+ t_k u^k$. By $o(t)$, we mean $\lim_{t\rightarrow0}\frac{o(t)}{t}=0$.
We denote by $\mathbb B$, $\bar{\mathbb B}$, $\mathbb S$  the open unit ball, the closed unit ball and 
  the unit sphere, respectively. $\mathbb B_\delta(\bar z)$ denotes the open unit ball centered at $\bar z$ with radius $\delta$.
We denote by ${\rm co}\Omega$ and $cl\Omega$ the convex hull and  the closure of a set $\Omega$, respectively. The distance from a point $x$ to a set $\Omega$ is denoted by  ${\rm dist}(x,\Omega):=\inf\{\|x-y\||y\in\Omega\}$ and the indicator function of set $\Omega$  is denoted by  $\delta_\Omega$. For a single-valued map $\phi:\mathbb R^n\rightarrow\mathbb R^m$, we denote by $\nabla \phi(x)\in \mathbb{R}^{m\times n}$   the Jacobian matrix of $\phi$  at $x$ and for a function $\varphi:\mathbb R^n\rightarrow\mathbb R$, we denote by $\nabla \phi(x)$ both the gradient and the Jacobian of $\phi$  at $x$. {Denote the pre-image of set $\Omega\subseteq\mathbb R^m$ under map $\phi$ by $\phi^{-1}(\Omega):=\{x\in\mathbb R^n|\phi(x)\in\Omega\}$.} For a set-valued map $\Phi:\mathbb R^n\rightrightarrows\mathbb R^m$ the graph of $\Phi$ is defined by  gph$ \Phi:=\{(x,y)| y\in \Phi(x)\}$. {For an extended-valued function $\varphi:\mathbb R^n\rightarrow \overline{\mathbb R}$, we define its domain by $dom\varphi:=\{x\in\mathbb R^n| \varphi(x) <\infty\}$, and its epigraph by $epi\varphi:=\{(x,\alpha)\in\mathbb R^{n+1}|\alpha\geq\varphi(x)\}$. For a function $g:\mathbb R^n\rightarrow\mathbb R$, we denote $g_+(x):=\max \{0,g(x)\}$ and if it is vector-valued then the maximum is taken componentwise.}


 We now review some basic concepts and results in variational analysis, which will be used later on. For more details see e.g. \cite{BHA,BS,Clarke,DonRock,Long,Aub2,RW}. Moreover we derive some preliminary results that will be needed.
\begin{defn}[Tangent Cone and Normal Cone] (see, e.g., \cite[Definitions 6.1 and 6.3]{RW})
	Given a set $\Omega\subseteq\mathbb{R}^n$ and a point $\bar{x}\in \Omega$, the tangent cone to $\Omega$ at $\bar{x}$ is defined as
	$$T_\Omega(\bar{x}):=\left \{d\in\mathbb{R}^n|\exists t_k\downarrow0, d_k\rightarrow d\ \mbox{ s.t. } \bar{x}+t_kd_k\in \Omega\ \forall k\right \}.$$
	The regular normal cone, the limiting normal cone  to $\Omega$ at $\bar{x}$ are defined as
\begin{eqnarray*}\widehat{N}_\Omega(\bar{x})&:=&\left \{ \zeta\in \mathbb{R}^n\bigg| \langle \zeta ,x-\bar{x}\rangle \leq o(\|x-\bar x\|) \quad \forall x\in \Omega \right \},\\
	 N_\Omega(\bar{x})&:=&\left \{\zeta\in \mathbb{R}^n\bigg| \exists \ x_k\xrightarrow{\Omega}\bar{x},\ \zeta_k{\rightarrow}\zeta\ \text{such that}\ \zeta_k\in\widehat{N}_\Omega(x_k)\ \forall k\right \},\end{eqnarray*}
	respectively.
	\end{defn}

	\begin{defn}[Directional Normal Cone] (\cite[Definition 2.3]{GM}  or \cite[Definition 2]{Gfr13}). 
    Given a set $\Omega \subseteq \mathbb{R}^n$, a point $\bar x \in \Omega$ and a direction $d\in\mathbb{R}^n$, the limiting normal cone to $\Omega$ at $\bar{x}$ in direction $d$ is defined by
    $$N_\Omega(\bar{x};d):=\left \{\zeta \in \mathbb{R}^n\bigg| \exists \ t_k\downarrow0, d_k\rightarrow d, \zeta_k\rightarrow\zeta  \mbox{ s.t. } \zeta_k\in \widehat{N}_\Omega(\bar{x}+t_kd_k)\ \forall k \right \}.$$
\end{defn}
It is obvious that $N_{\Omega}(\bar x; 0)=N_{\Omega}(\bar x)$, $N_{\Omega}(\bar x; d)=\emptyset$ if $d \not \in T_\Omega(\bar x)$ and $N_{\Omega}(\bar x;d)\subseteq N_\Omega(\bar x)$.
Moreover when $\Omega$ is convex, by \cite[Lemma 2.1]{Gfr14} the directional and the classical normal cone have the following relationship
\begin{equation}
N_\Omega(\bar x;d)=N_\Omega(\bar x)\cap \{d\}^\perp  \qquad \forall d\in T_\Omega(\bar x).\label{convNormal}
\end{equation}


When $u=0$ the following  definition coincides with the Painlev\'e-Kuratowski inner/lower and outer/upper  limit of $\Phi$ as $x\rightarrow\bar x$ respectively; see e.g., \cite{Aub2}.
\begin{defn}\label{innerlimit} Given a set-valued map $\Phi:\mathbb R^n\rightrightarrows\mathbb R^m$ and a direction $
d \in\mathbb R^n$, the inner/lower and outer/upper limit of $\Phi$ as $x\xrightarrow{d}\bar x$ respectively is defined by
	\begin{align*}
	\liminf_{x\xrightarrow{d}\bar x} \Phi(x):=\{y\in\mathbb R^m|&\forall \ \mbox{sequence}\ t_k\downarrow0, d^k\rightarrow d, \exists y^k\rightarrow y \mbox{ s.t. } y^k\in \Phi(\bar x+t_kd^k)\}\\
	\limsup_{x\xrightarrow{d}\bar x} \Phi(x):=\{y\in\mathbb R^m|&\exists   \ \mbox{sequence}\ t_k\downarrow0, d^k\rightarrow d,  y^k\rightarrow y \mbox{ s.t. } y^k\in \Phi(\bar x+t_kd^k)\},
	\end{align*}
	respectively.
\end{defn}

\begin{defn}[Directional derivatives]
Let $\phi:\mathbb R^n\rightarrow\mathbb{R}^m $ and $x, u\in \mathbb R^n$. The usual directional derivative of $\phi$ at $x$ in the direction $u$ is 
$$\phi'(x;u):=\lim_{t\downarrow0}\frac{\phi(x+tu)-\phi(x)}{t}$$
when this limit exists. 
\end{defn}
We say that $\phi:\mathbb R^n\rightarrow\mathbb{R}^m $ is directionally Lipschitz continuous  at $\bar x$ in direction $u$ if there are positive numbers $L, \epsilon, \delta$ such that
$$\|\phi(x)-\phi(x') \|\leq L \|x-x'\|  \quad \forall x, x' \in \bar x+{\cal V}_{\epsilon,\delta}(u). $$
It is easy to see that if $\phi:\mathbb R^n\rightarrow\mathbb{R}^m $ is directionally Lipschitz continuous and directionally differentiable at $x$ in direction $u$
then for all sequence $\{u^k\}$ which converges to $u$, we have
\begin{equation*}
\phi'(x;u)=\lim_{k\rightarrow\infty}\frac{\phi(x+t_ku^k)-\phi(x)}{t_k}.
\end{equation*}

We now recall the definition of some subdifferentials below.
\begin{defn}[Subdifferentials](\cite[Definition 8.3]{RW})
Let  $\varphi:\mathbb{R}^n \rightarrow \overline{\mathbb{R}}$ and  $\bar x\in  {\rm dom} \varphi$. 
The Fr\'{e}chet (regular) subdifferential of $\varphi$ at $\bar{x}$ is the set
\begin{eqnarray*}
\widehat\partial \varphi(\bar{x}):=\left\{\xi\in\mathbb{R}^n|\varphi(x)\geq \varphi(\bar x)+\langle\xi,x-\bar x\rangle+{o(\|x-\bar x\|)}\right\},
\end{eqnarray*}
the limiting (Mordukhovich or basic) subdifferential of $\phi$ at $\bar{x}$ is the set
\begin{align*}
\partial \varphi(\bar{x}):=
\{\xi\in\mathbb{R}^n|\exists x^k\rightarrow  \bar x,  \ \xi^k\rightarrow\xi \ \mbox{ s.t. }  \varphi(x^k)\rightarrow \varphi(\bar x), \xi^k\in \widehat\partial \varphi(x^k)\}.
\end{align*}
\end{defn}
\begin{defn}[Directional subdifferentials]\label{ads}
	 \cite{Gfr13,GM,Long,BHA} Let $\varphi:\mathbb{R}^n\rightarrow \overline{\mathbb{R}}$ and $\bar x\in dom\varphi$. 
The limiting subdifferential of $\varphi$ at $\bar x$ in direction 
$u\in \mathbb{R}^n$ is defined as 
	\begin{eqnarray*}
	\partial \varphi(\bar x;u):=\{\xi\in\mathbb{R}^n|\exists t_k\downarrow0, u^k\rightarrow u, \xi^k\rightarrow\xi \mbox{ s.t. } \varphi(\bar x+t_ku^k)\rightarrow \varphi(\bar x), \   \xi^k\in \widehat{\partial} \varphi(\bar x+t_ku^k)\}.
	\end{eqnarray*} 
\end{defn} It is easy to see that if $u\notin T_{dom \varphi}(\bar x)$, then $ \partial \varphi(\bar x;u)=\emptyset$ and $\partial \varphi(\bar x;0)=\partial \varphi(\bar x)$.

\begin{defn}[Directional Clarke subdifferential]\label{Clarke}
Let $\varphi:\mathbb R^n\rightarrow\mathbb R$ be directionally Lipschitz continuous at  $\bar x$ in direction $u\in\mathbb R^n$. The directional Clarke subdifferential of $\varphi$ at $\bar x$ in direction $u$ is defined as
\begin{equation*}
\partial^c \varphi(\bar x;u):=co(\partial \varphi(\bar x;u)).
\end{equation*}	
\end{defn}
 It is clear that  the directional Clarke subdifferential in direction $u=0$ coincides with the Clarke   subdifferential $\partial^c \varphi(\bar x)$.
\begin{prop}\label{Prop2.4} Let $\varphi:\mathbb{R}^n\rightarrow \mathbb{R}$ be directionally Lipschitz continuous at $\bar x$ in direction $u$. Then we have
$$\partial^c \varphi(\bar x;u)=co\limsup_{x\xrightarrow{u}\bar x}\partial^c \varphi(x),\ \partial^c(-\varphi)(\bar x;u)=-\partial^c \varphi(\bar x;u).$$
\end{prop}
\beginproof By \cite[Theorem 5.4]{Long},  we have
$\partial \varphi(\bar x;u)=\displaystyle \limsup_{x\xrightarrow{u}\bar x} \partial \varphi(x).$
{It follows that 
$$\partial^c \varphi(\bar x;u):=co(\partial \varphi(\bar x;u)) = co (\limsup_{x\xrightarrow{u}\bar x} \partial \varphi(x)).$$
Hence  to prove the first equality stated in the proposition, we only need to show that 
\begin{equation}\label{limsup}
   \limsup_{x\xrightarrow{u}\bar x} co(\partial \varphi(x)) \subseteq co (\limsup_{x\xrightarrow{u}\bar x} \partial \varphi(x)).
\end{equation}
Let  $\zeta\in \displaystyle \limsup_{x\xrightarrow{u}\bar x} co( \varphi(x)) $. Then there exist  sequences $x^k\xrightarrow{u}\bar x$ and  $\zeta^k\rightarrow \zeta$ such that $\zeta^k\in co( \varphi(x^k))$. By the Carath\'eodory Theorem, for each $k$, there exist $\{a^k_1,a^k_2,\ldots,a^k_{n+1}\}\subseteq \partial\varphi(x^k)$ and nonnegative scalars $\{\lambda^k_1,\lambda^k_2,\ldots,\lambda^k_{n+1}\}$ such that $\zeta^k=\Sigma_{i=1}^{n+1}\lambda^k_ia^k_i$ and $\Sigma_{i=1}^{n+1}\lambda^k_i=1$. Since $\varphi(x)$ is directionally Lipschitz continuous at $\bar x$ in direction $u$, by \cite[Theorem 9.13]{RW}, there exists $L>0$ such that $\partial \varphi(x^k)\subseteq L\bar{\mathbb B}$ for all $k$ sufficiently large. Hence the sequence $\{(a^k_1,a^k_2,\ldots,a^k_{n+1},\lambda^k_1,\lambda^k_2,\ldots,\lambda^k_{n+1})\}$ is bounded. Passing to a subsequence if necessary, we have  $(a^k_1,a^k_2,\ldots,a^k_{n+1},\lambda^k_1,\lambda^k_2,\ldots,\lambda^k_{n+1})\rightarrow(a_1,a_2,\ldots,a_{n+1},\lambda_1,\lambda_2,\ldots,\lambda_{n+1})$ as $k\rightarrow\infty$. Consequently, we have for each $i=1,2,\ldots,n$, $\lambda_i\geq0$, $a_i\in \displaystyle \limsup_{x\xrightarrow{u}\bar x}\varphi(x)$, $\Sigma_{i=1}^{n+1}\lambda_i=1$ and $\Sigma_{i=1}^{n+1}\lambda_ia_i=\zeta$. This implies that $\zeta\in co \displaystyle(\limsup_{x\xrightarrow{u}\bar x} \partial \varphi(x))$. Hence (\ref{limsup}) holds.
The second equality stated in the proposition follows directly from the first equality and the scalar multiplication rule of Clarke subdifferential \cite[Proposition 2.3.1]{Clarke}.} 
\endproof

We now give the definition of directional metric subregularity constraint qualification.
\begin{defn}[Directional MSCQ] \cite[Definition 2.1]{Gfr13}\label{dims}
Let $\bar x$ be a solution to the system $g(x)\leq0$, where $g:\mathbb R^n\rightarrow\mathbb R^m$.
 Given a direction $u\in\mathbb{R}^n$, the system $g(x)\leq0$ is said to satisfy the directional metric subregular constraint qualification (MSCQ) at $\bar{x}$ in direction $u$, if  there are positive reals $\epsilon>0,\delta>0,$ and $\kappa>0$ such that
		\begin{equation*}
		{{\rm dist}}(x,g^{-1}(\mathbb R^m_-))\leq\kappa \|g_+(x)\| \qquad \quad \forall  x\in\bar{x}+{\cal V}_{\epsilon,\delta}(u).
		\end{equation*}
\end{defn}
If $u=0$ in the above definition, then we say that the system $g(x)\leq0$ satisfies MSCQ at $\bar x$. 

\section{Directional KKT conditions under directional  calmness condition}\label{Sec3}
In this section we derive directional KKT condition for  the optimization problem
\begin{eqnarray*}
{\rm (P)}~~~ \min_{z} && \varphi(z)\quad \mbox{ s.t. } \phi(z)\leq 0,
 \end{eqnarray*}
{where $\varphi:\mathbb{R}^n \rightarrow \mathbb{R} $ and $\phi:\mathbb{R}^n \rightarrow \mathbb{R}^q $.}

{The concept of (Clarke) calmness for a mathematical program was first defined by Clarke \cite[Definition 6.41]{Clarke}.  We now introduce a directional version of the calmness condition for (P).}
	\begin{defn}[Directional Clarke calmness]
		Suppose $\bar z$ solves {\rm (P)}. We say that {\rm (P)} is (Clarke) calm at $\bar z$ in direction $u$ if there exist positive scalars $\epsilon,\delta,\rho$, such that for any $\alpha \in\epsilon\mathbb B$ and any $z\in\bar z+{\cal V}_{\epsilon,\delta}(u)$ satisfying $\phi(z)+\alpha\leq 0$ one has,
		\begin{equation*}
		\varphi(z)-\varphi(\bar z)+\rho\|\alpha\|\geq0.
		\end{equation*}
	\end{defn}
We now prove that the directional MSCQ implies the directional calmness of problem (P) provided the objective function is directionally Lipschitz continuous.
	\begin{lemma}\label{mscalm}
		Let $\bar z$ solve $({\rm P})$ and $\varphi(z)$ be directionally Lipschitz continuous at $\bar z$ in direction $u$. Suppose  that the system $\phi(z)\leq0$ satisfies the directional MSCQ at $\bar z$ in direction $u$. Then $({\rm P})$ is calm at $\bar z$ in direction $u$.
	\end{lemma}
\beginproof
Since $\phi(z)\leq0$ satisfies the directional MSCQ at $\bar z$ in direction $u$, by Definition \ref{dims}, there exist positive scalars $\epsilon,\delta,\kappa$ such that
\begin{equation}\label{2}
{\rm dist}(z,\phi^{-1}(\mathbb R^q_-))\leq\kappa\|\phi_+(z)\|\qquad \forall z\in\bar z+\mathcal V_{\epsilon,\delta}(u).
\end{equation} Let $\tilde z$ be the projection of $z$ onto $\phi^{-1}(\mathbb R^q_-)$. Since $\varphi(z)$ is directionally Lipschitz continuous, without loss of generality, taking $\epsilon,\delta$ small enough, there exists $L>0$ such that $|\varphi(z)-\varphi(z')|\leq L\|z-z'\|$ for any $z,z'\in\bar z+\mathcal V_{\epsilon,\delta}(u)$. Then we have for any  $\alpha\in\epsilon\mathbb B$ satisfying $\phi(z)+\alpha\leq 0$,
\begin{align*}
\varphi(z)-\varphi(\bar z)+L\kappa\|\alpha\|
&\geq\varphi(z)-\varphi(\bar z)+L\kappa\|\phi_+(z)\|\\
&\geq\varphi(z)-\varphi(\bar z)+L{\rm dist}(z,\phi^{-1}(\mathbb R^q_-))\\
&=\varphi(z)-\varphi(\bar z)+L\|z-\tilde z\|\\
&\geq\varphi(z)-\varphi(\tilde z)+L\|z-\tilde z\|\\
&\geq0,
\end{align*} where the second inequality follows from (\ref{2}), the third inequality follows from the optimality of $\varphi(z)$ at $\bar z$ and the last inequality follows from the directional Lipschitz continuity of $\varphi(z)$ at $\bar z$. Let $\rho:=L\kappa$. The proof is complete.
\endproof

Let $\bar z$ be a feasible solution to problem (P). We denote by $$\bar{I}_\phi:=I_\phi(\bar z):=\{j=1,\ldots,q|\phi_j(\bar z)=0\}$$ the set of indexes of active constraints at $\bar z$. If $\varphi$ is continuously differentiable and $\phi$ is directionally Lipschitz continuous and directionally differentiable,  we define the linearized cone by $L(\bar z):=\{u\in\mathbb R^n|\phi'_j(\bar z;u)\leq0, j\in I_\phi(\bar z)\}$ and  the critical cone by $$C(z):=\{u\in L(\bar z)| \nabla\varphi(z)u\leq0\}=\{u\in\mathbb R^n|\phi'_j(\bar z;u)\leq0, j\in I_\phi(\bar z), \nabla\varphi(z)u\leq0\}.$$

{The following definition lists  some sufficient conditions for the directional MSCQ, hence are sufficient for directional calmness.}
\begin{defn}\label{qp} Let $\phi (\bar z) \leq 0$ and $u\in \mathbb{R}^n$.
\begin{itemize}
		\item Suppose that  $\phi$ is Lipschitz at $\bar z$. We say that  the no-nonzero abnormal multiplier constraint qualification (NNAMCQ) holds at $\bar z$ if 
\begin{equation*}
0\in \partial\langle\zeta,\phi\rangle(\bar z) \mbox{ and } 0\leq\zeta\perp \phi(\bar z)  \Longrightarrow \zeta=0.
\end{equation*}
\item  Suppose that   $ \phi$ is directionally Lipschitz and  directionally differentiable at $\bar z$ in direction $u$.
  We say that the first order sufficient condition for metric subregularity (FOSCMS)  holds  at $\bar z$ in direction $u$ if there exists no $\zeta
\neq 0$ satisfying $0\leq\zeta\perp \phi(\bar z),\ \zeta\perp \phi'(\bar z;u)$ and
	\begin{equation}\label{FOSCMScon}
	0\in \partial\langle\zeta,\phi\rangle(\bar z;u).
	\end{equation}	
	\item Suppose that  $\phi$ is directionally Lipschitz and  directionally differentiable at $\bar z$ in direction $u$.
We say that the directional quasi-normality  holds at $\bar{z}$ in direction $u$
if 
 there exists no $\zeta\neq 0$ satisfying $0\leq\zeta\perp \phi(\bar z),\ \zeta\perp \phi'(\bar z;u)$ such that
 (\ref{FOSCMScon}) holds 
		and there exists sequences $t_k\downarrow0,\ u^k\rightarrow u$ satisfying
		\begin{equation}
 \phi_j(\bar{z}+t_ku_k)>0,\ \mbox{ if }  j \in \bar{I}_\phi  \mbox{ and } \zeta_j>0\label{qua 1}.
		\end{equation}
		\end{itemize}
\end{defn}
 
\begin{prop}\label{qpms}  Let $\phi(\bar z)\leq0$ and suppose that $ \phi$ is directionally Lipschitz and  directionally differentiable at $\bar z$ in direction {$u\in L(\bar z)$.}
	If the directional quasi-normality holds at $\bar z$ in direction $u$ for the inequality system $\phi(z)\leq 0$. Then  the system $\phi(z)\leq0$ satisfies the directional MSCQ at $\bar z$ in direction $u$.
\end{prop}
\beginproof  Since $\phi$ is directionally Lipschitz  and directionally differentiable at $\bar z$ in direction $u$,
by \cite[Corollary 4.1, Proposition 5.1]{BHA} we have  $$D^*\phi(\bar z;(u,\phi'(\bar z;u)))(\zeta)=\partial\langle\zeta,\phi\rangle(\bar z;u),$$
where $D^*\phi(\bar z;(u,v))$ is the limiting coderivative of $\phi$ at $\bar z$ in direction direction $(u,v)$ as defined in \cite{BHA}.
By equality $(\ref{convNormal})$, we have $N_{\mathbb R^q_-}(\phi(\bar z);\phi'(\bar z;u))=\{\mu\in\mathbb R^q|0\leq\mu\perp \phi(\bar z), \mu\perp \phi'(\bar z;u)\}$. {For any sequences $\{\zeta^k\}$ and $\{s^k\}$}, if $\zeta_j>0$ and $\widehat{N}_{\mathbb R^q_-} (s^k_j) \ni \zeta^k_j \rightarrow \zeta_j$, then for large enough $k$,  $\zeta^k_j>0$ and hence $s^k_j=0$. Hence the condition $(\ref{qua 1})$ is equivalent to the sequential condition in \cite[Definition 4.1(a)]{BYZ}. 
Therefore the quasi-normality in direction $u\in L(\bar z)$ means that  there exists no $\zeta \not =0$ such that
$$0\in D^* \phi (\bar z; (u,\phi'(\bar z;u)))(\zeta), \qquad \zeta \in N_{\mathbb R^q_-}(\phi(\bar z);\phi'(\bar z;u))$$
and there exist sequences $t_k\downarrow0,\ u^k\rightarrow u$ such that (\ref{qua 1}) holds.

 From the proof of \cite[Lemma 3.1 and Corollary 4.1]{BYZ} and \cite[Corollary 1]{Hcor}, one can easily obtain that the quasi-normality at $\bar z$ in direction $u$ implies that $\phi(z)\leq0$ satisfies the MSCQ at $\bar z$ in direction $u$.
\endproof

{In the following theorem, we derive the directional KKT condition under the directional  calmness condition.}
\begin{thm}\label{dKKT}
	Let $\bar{z}$ be a local minimizer of {\rm (P)}. Suppose that $\varphi(z)$ is continuously differentiable at $\bar z$ and $\phi(z)$ is directionally Lipschitz and directionally differentiable at $\bar z$ in direction $u\in 
	C(\bar z).$ 
	{Suppose that problem ${\rm (P)}$ is calm at $\bar z$ in direction $u$.} Then there exists a vector $\lambda_\phi\in\mathbb R^{q}$ such that
	$ 0\leq\lambda_\phi\perp \phi(\bar z),\ \lambda_\phi\perp \phi'(\bar z;u)$ and 
	\begin{eqnarray*}
	0\in\nabla \varphi(\bar z)+\partial\langle\lambda_\phi,\phi\rangle(\bar z;u).
	\end{eqnarray*} 
\end{thm}
\beginproof
Since (P) is calm at $\bar z$ in direction $u$, there exist positive scalars $\epsilon,\delta,\rho$ such that
\begin{equation}\label{penalizedP}
\varphi(z)+\rho\|\phi_+(z)\|\geq \varphi(\bar z)\qquad  \forall z\in \bar z+{\cal V}_{2\epsilon,\delta}(u).
\end{equation} 
Since $u\in C(\bar z)$, for $t>0$ sufficiently small,  we have  $\phi(\bar z)+t\phi'(\bar z;u)\leq0$ and hence
$$0\leq \frac{\|\phi_+(\bar z+tu)\|}{t}\leq \frac{\|\phi(\bar z+tu)-\phi(\bar z)-t\phi'(\bar z;u)\|}{t} .$$
It follows that $\lim_{t\downarrow 0} \frac{\|\phi_+(\bar z+tu)\|}{t}=0$.

Since  $\bar z+t u\in \bar z+cl({\cal V}_{\epsilon,\delta}(u))$ for $t$ sufficiently small,   by (\ref{penalizedP}),  $$\varphi(\bar z+tu)+\rho\|\phi_+(\bar z+tu)\|\geq\varphi(\bar z)$$ for all $t$ small enough. Together with 
$\nabla \varphi(\bar z)u\leq0$ we have
\begin{equation}\label{speed}
\lim_{t\downarrow0}\frac{\varphi(\bar z+tu)+\rho\|\phi_+(\bar z+tu)\|-\varphi(\bar z)}{t}=0.
\end{equation}
For each $k=0,1,\ldots$, define $\sigma_k:=2(\varphi(\bar z+\frac{u}{k})+\rho\|\phi_+(\bar z+\frac{u}{k})\|-\varphi(\bar z))$. 
If $\sigma_k\equiv0$, then for each large enough $k$, by (\ref{penalizedP}), $\bar z+\frac{u}{k}$ is a global minimizer of the function $\varphi(z)+\rho\|\phi_+(z)\|+\delta_{\bar z+cl({\cal V}_{\epsilon,\delta}(u))}(z)$. Since for each large enough $k$, $\bar z+\frac{u}{k}$ is an interior point of $\bar z+cl({\cal V}_{\epsilon,\delta}(u))$, by the well-known Fermat's rule and the calculus rule ({see e.g., \cite[Corollary 10.9]{RW}}), 
\begin{equation}\label{opk1}
0\in\nabla\varphi(\bar z+\frac{u}{k})+\rho\partial(\|\phi_+{(\cdot)}\|)(\bar z+\frac{u}{k}).
\end{equation}
Otherwise, without loss of generality, we assume that  for all $k$, $\sigma_k>0$.  Then by definition of $\sigma_k$ we have for $k$ sufficiently large, $$\varphi(\bar z+\frac{u}{k})+\rho\|\phi_+(\bar z+\frac{u}{k})\|+\delta_{\bar z+cl({\cal V}_{\epsilon,\delta}(u))}(\bar z+\frac{u}{k})<\varphi(\bar z)+\sigma_k.$$  Define $\lambda_k:=\frac{2\|u\|r}{k\epsilon}\sqrt{\frac{\sigma_kk\epsilon}{2\|u\|r}}$. By Ekeland's variation principle (see e.g., \cite[Theorem 2.26]{Aub2}), there exists $\tilde z^k$ satisfying that $\|\tilde z^k-(\bar z+\frac{u}{k})\|\leq\lambda_k$, and the function
$$z\rightarrow \varphi(z)+\rho\|\phi_+( z)\|+\delta_{\bar z+cl({\cal V}_{\epsilon,\delta}(u))}(z)+\frac{\sigma_k}{\lambda_k}\|z-(\bar z+\frac{u}{k})\|$$ attains its global minimum at $\tilde z^k$. Since $\frac{\epsilon u}{2\|u\|}$ is an interior point of $cl({\cal V}_{\epsilon,\delta}(u))$,  there exists $r\in (0, \epsilon/2)$ such that $\frac{\epsilon u}{2\|u\|}+r\mathbb{B} \subseteq cl({\cal V}_{\epsilon,\delta}(u))$. It is obvious that the following implication holds
{\begin{equation*}
z\in cl({\cal V}_{\epsilon,\delta}(u)), 0\leq\alpha\leq1 \Rightarrow\alpha z \in  cl({\cal V}_{\alpha\epsilon,\delta}(u))
\end{equation*}} Hence  $(\frac{\epsilon u}{2\|u\|}+r\mathbb{B} ) \frac{2\|u\|}{\epsilon k}\subseteq cl({\cal V}_{\epsilon,\delta}(u))$ and hence  $\bar z+\frac{u}{k}+\frac{2\|u\|}{k\epsilon}r\mathbb B\subseteq (\bar z+cl({\cal V}_{\epsilon,\delta}(u)))$ and since $\sigma_k=o(\frac{1}{k})$  by (\ref{speed}), $\tilde z^k$ is in the interior of $\bar z+cl({\cal V}_{\epsilon,\delta}(u))$. Then by the well-known Fermat's rule, we obtain
\begin{equation}\label{opk2}
0\in\nabla\varphi(\tilde z^k)+\rho\partial(\|\phi_+(\cdot)\|)(\tilde z^k)+\frac{\sigma_k}{\lambda_k}\bar{\mathbb B}.
\end{equation}
Since $\phi$ is directionally Lipschitz continuous at $\bar z$ in direction $u$, it is Lipschitz continuous at $\tilde z^k$ for $k$ large enough. So  by the chain rule for limiting subdifferential \cite[Corollary 3.43]{Aub2}, we have  $$\partial(\|\phi_+(\cdot)\|)(\tilde z^k)\subseteq\cup_{\zeta'\in  \partial\|(\cdot)_+\|(\phi(\tilde z^k))}\partial\langle\zeta',\phi\rangle(\tilde z^k).$$
Therefore by  $(\ref{opk1})$ or $(\ref{opk2})$, $\exists\zeta^k\in\partial\| (\cdot)_+\|(\phi(\bar z+\frac{u}{k}))$ or $\exists\zeta^k\in\partial\| (\cdot)_+\|(\phi(\tilde z^k))$ such that
\begin{equation}\label{opk}
0\in\nabla\varphi(\bar z+\frac{u}{k})+\rho\partial\langle\zeta^k,\phi\rangle(\bar z+\frac{u}{k}),\ \mbox{or}\
0\in\nabla\varphi(\tilde z^k)+\rho\partial\langle\zeta^k,\phi\rangle(\tilde z^k)+\frac{\sigma_k}{\lambda_k}\bar{\mathbb B}.
\end{equation}
Since the function $\|x_+\|$ is  Lipschitz continuous, by \cite[Theorem 9.13]{RW}, $\{\zeta^k\}$ is bounded. Without loss of generality, there exists $\zeta:=\lim_k\zeta^k$. By the way, one can easily obtain that $\lim_k(\bar z+u/k-\bar z)/\frac{1}{k}=\lim_k(\tilde z^k-\bar z)/\frac{1}{k}=u$. Since $\sigma_k=o(\frac{1}{k}),\ \lim_k\frac{\sigma_k}{\lambda_k}=0$. Taking the limit of $(\ref{opk})$ as $k\rightarrow\infty$,   by \cite[Theorem 5.4]{Long} we have
$$0\in \nabla \varphi(\bar z)+\rho \partial \langle \zeta, \phi \rangle (\bar z; u).$$
Moreover by  \cite[Corollary 4.2]{BHA}, $\zeta\in\partial(\|(\cdot)_+\|)(\phi(\bar z);\phi'(\bar z;u))\subseteq N_{\mathbb R^q_-}(\phi(\bar z);\phi'(\bar z;u))$.
The desired result holds by taking $\lambda_{\phi}:=\rho\zeta\in N_{\mathbb R^q_-}(\phi(\bar z);\phi'(\bar z;u))=\{\xi\in\mathbb R^q|0\leq\xi\perp\phi(\bar z), \xi\perp\phi'(\bar z;u)\}$. 
\endproof
{We now give an example of a bilevel program where  the partial calmness and calmness fail but the   calmness condition holds in a nonzero critical direction.}
\begin{example}\label{Ex3.1} Consider the following bilevel program:
		\begin{eqnarray*}
		\mbox{(BP)}\quad \min && F(x,y):=(x-y-1)^{\frac{5}{3}}+4(x+y+1)^{\frac{5}{3}}\\
		s.t. &&\ -1\leq x\leq 1, y\in S(x),
\end{eqnarray*}
where for each $x$, $S(x)$ is the solution set for the lower level program:
$$   \min_{y} \{f(x,y):=-(x+y)^2+x^3(x+y-1),\\
	 s.t.\ -y-x-1\leq0,
	 y+x-1\leq0\}.$$
It is easy to see that the solution mapping $S(x)$ of the lower level problem is equal to
\begin{eqnarray}\label{solutions}
S(x)=\left\{\begin{array}{ll}
-x-1,\ &x>0,\\
\{-1,1\},\ &x=0,\\
-x+1,\ &x<0.
\end{array}
\right.
\end{eqnarray}  And the global optimal solution of (BP) is $(\bar x,\bar y)=(0,-1)$. The constraints 
$ y+x-1\leq0$ and $-1\leq x\leq1$ are inactive at $(0,-1)$. The value function 
\begin{equation}
V(x)=\left \{ \begin{array}{ll}
-1-2x^3 & x>0,\\
-1 & x\leq 0.
\end{array}
\right.
\label{valuefunction}
\end{equation}

First, we prove that the partial calmness condition fails at $(\bar x,\bar y)$. For any scalar $\rho>0$, consider the partially penalized problem:
\begin{eqnarray*}
	(VP)_\rho\ &\min& F(x,y)+\rho(f(x,y)-V(x))\\
	&s.t.&\ g_1(x,y):=-y-x-1\leq0,
	g_2(x,y):= y+x-1\leq0,\\
	&&-1\leq x\leq 1.
\end{eqnarray*}
Since $-1<\bar x<1, g_1(\bar x,\bar y)=0,\ g_2(\bar x,\bar y)<0$, by  (\ref{solutions})-(\ref{valuefunction}),  the critical cone is
\begin{align*}
C(\bar x,\bar y)&=\{(u,v)|\nabla F(\bar x,\bar y)(u,v)\leq0, \nabla f(\bar x,\bar y)(u,v)-V'(\bar x;u)=0,\nabla g_1(\bar x,\bar y)(u,v)\leq 0\}\\
&= \{(u,v)|u+v=0\}.
\end{align*}
Consider the sequence  $(x^k,y^k):=(-\frac{1}{k},\frac{1}{k}-1)$ which are feasible to $(VP)_\rho$ and converges to $(\bar x,\bar y)$.  Since $F(x^k,y^k)=-(\frac{2}{k})^{\frac{5}{3}},\ f(x^k,y^k)=-1+\frac{2}{k^3}$ and by $(\ref{valuefunction})$, $V(x^k)=-1$,  we have $F(x^k,y^k)+\rho(f(x^k,y^k)-V(x^k))=-(\frac{2}{k})^{\frac{5}{3}}+\frac{2\rho}{k^3}$. Hence  for $k$ sufficiently large, we have $$F(x^k,y^k)+\rho(f(x^k,y^k)-V(x^k))<0=F(\bar x,\bar y)+\rho(f(\bar x,\bar y)-V(\bar x)).$$ This means that for any $\rho>0$, $(\bar x,\bar y)$ is not a local minimizer of $(VP)_\rho$. Hence, the partial calmness fails. Since the calmness condition is in general stronger than partial calmness, the calmness condition also fails. In fact for this example since the constraint functions for $(VP)_\rho$ are all affine, the partial calmness is equivalent to the fully calmness.
Notice that $(x^k,y^k)\rightarrow (\bar x,\bar y)$ in direction $(-1,1)$ and so we have shown that problem (VP) is not calm in direction $(-1,1)$.
{Next, we prove that \ (VP) is calm at $(\bar x,\bar y)$ in direction $(1,-1)\in C(\bar x,\bar y)$.}
Since the constraints $g_2(x,y)\leq 0$ and $-1\leq x\leq 1$ are inactive at $(\bar x,\bar y)=(0,-1)$, it suffices to show that
 there exists a positive scalar $\rho$ such that for any sequences $t_k\downarrow0,\ (u^k,v^k)\rightarrow (\bar u,\bar v):=(1,-1)$, 
for $k$ sufficiently large, 
\begin{align}
F(\bar x+t_ku^k,\bar y+t_kv^k)+\rho\|(&f(\bar x+t_ku^k,\bar y+t_kv^k)-V(\bar x+t_ku^k))_+\|\notag\\
&+\rho\|(g_1)_+(\bar x+t_ku^k,\bar y+t_kv^k)\|-F(\bar x,\bar y)\geq 0.\label{dcal}
\end{align}
Suppose that  $g_1(\bar x+t_ku^k,\bar y+t_kv^k)\leq 0$. Then for $k$ sufficiently large, $\bar y+t_kv^k$ is a feasible solution for $(P_{\bar x+t_ku^k})$ and hence  $f(\bar x+t_ku^k,\bar y+t_kv^k)-V(\bar x+t_ku^k)\geq0$ by the definition of the value function. Moreover $F(\bar x+t_ku^k,\bar y+t_kv^k)\geq F(\bar x,\bar y)$. Hence (\ref{dcal}) holds. Otherwise suppose that  $$g_1(\bar x+t_ku^k,\bar y+t_kv^k)=-t_k(u^k+v^k)>0.$$
Hence $t_k(u^k+v^k)<0$. Together with $u^k>0, t_k>0$,  we can verify that   $f(\bar x+t_ku^k,\bar y+t_kv^k)-V(\bar x+t_ku^k)<0.$ Also since  $t_k\downarrow0, -(u^k+v^k)\downarrow0$, we have \begin{align*}
&F(\bar x+t_ku^k,\bar y+t_kv^k)-F(\bar x,\bar y)+\rho\|(g_1)_+(\bar x+t_ku^k,\bar y+t_kv^k)\|\\
=&t_k^{\frac{5}{3}}(u^k+v^k)^{\frac{5}{3}}-t_k(u^k+v^k)\geq0.
\end{align*}
 Hence, we obtain $(\ref{dcal})$.
%
Consequently {\rm (VP)} is calm at $(\bar x,\bar y)$ in direction $(\bar u,\bar v)=(1,-1)$.
\end{example}

\section{Directional sensitivity analysis of the value function}
In this section we study the directional sensitivity analysis of the value function of the lower level program $(P_x)$. The results of this section could be of independent interest.
Denote  the feasible map of the problem $(P_x)$ by 
$$\mathcal F(x):=\{y\in\mathbb R^m|g(x,y)\leq0\}$$ and the active index set $I_g(x,y):=\{i=1,\ldots,p|g_i(x,y)=0\}$. Let  the Lagrange function of $(P_x)$ be 
\[
{\cal L}(x,y;\lambda):=f(x,y)+g(x,y)^T\lambda
\]
and  the set of Lagrange multipliers associated with $y\in{\cal F}(x)$ be
\[
\Lambda(x,y):=\{\lambda\in\mathbb R^p|\nabla_y{\cal L}(x,y;\lambda)=0,\ g(x,y)^T\lambda=0,\ \lambda\geq0\}.
\]

\subsection{Preliminary results}
In this subsection we give some preliminary results that  may be needed. 
First we introduce a directional version of the restricted inf-compactness condition which was first introduced in \cite[Hypothesis 6.5.1]{Clarke} with the terminology  introduced in \cite[Definition 3.8]{GLYZ}. 
\begin{defn}[Directional Restricted Inf-compactness]
	We say that the restricted inf-compactness holds at $\bar x$ in direction $u$ with compact set $\Omega_u \subseteq \mathbb{R}^n$ if $V(\bar x)$ is finite and there exist  positive numbers $\epsilon>0, \delta>0$ such that for all $ x\in\bar x+{\cal V}_{\epsilon,\delta}(u)$ with $V(x)<V(\bar x)+\epsilon$, one always has $S(x)\cap\Omega_u\neq\emptyset$. When $u=0$ in the above, we say that the restricted inf-compactness holds at $\bar x$.
\end{defn}
Next we introduce a directional version of the inf-compactness condition  (see e.g., \cite[Page 272]{BS}).
It is not difficult to verify that the directional inf-compactness implies the directional restricted inf-compactness.
\begin{defn}[Directional Inf-compactness]
	We say that the inf-compactness holds at $\bar x$ in direction $u$ if there exist a compact set  $\Lambda_u\subseteq \mathbb{R}^n$, and positive numbers $\alpha>V(\bar x),\ \epsilon,\ \delta$ such that for all $x \in\bar x+{\cal V}_{\epsilon,\delta}(u)$, one always has {$\emptyset \neq \{y|f(x,y)\leq\alpha,\ g(x,y)\leq0 \}\subseteq \Lambda_u$}. When $u=0$ in the above, we say that the  inf-compactness holds at $\bar x$.
\end{defn}

When $f(x,y)$ satisfies the growth condition at $\bar x\in {\rm dom}\mathcal F$, i.e., there exists $\delta>0$ such that the set
\[
\{y\in\mathbb R^m| g(\bar x,y)\leq\alpha,\ f(\bar x,y)\leq M,\ \alpha\in\delta\bar{\mathbb B}\}
\]
is bounded for each $M\in\mathbb R$, the inf-compactness holds at $\bar x$. Similarly, if $f(x,y)$ is coercive or level bounded, the inf-compactness holds at $\bar x$.

{
\begin{defn}[Directional Inner Semi-compactness]\label{ic}
	Suppose $S(\bar x)\neq\emptyset$. We say that $S(x)$ is inner semi-compact at $\bar x$ in direction $u$ if for any sequences $t_k\downarrow0, u^k\rightarrow u$, there exists a sequence $y^k\in S(\bar x+t_ku^k)$ such that $\{y^k\}$ contains a convergent subsequence as $k\rightarrow\infty$.
	 When $u=0$ in the above, we say that $S(x)$ is inner semi-compact at $\bar x$.
\end{defn}
If $S(x)$ is inner semi-compact at $\bar x$ in direction $u$, then the inf-compactness holds at $\bar x$ in direction $u$. This is because otherwise, either (i) for $k\rightarrow\infty$, there exist sequences $t_k\downarrow0$ and $u^k\rightarrow u$ such that for any $y^k\in S(\bar x+t_ku^k)$, $f(\bar x+t_ku^k,y^k)\geq k$, or (ii) for $k\rightarrow\infty$, there exist sequences $t_k\downarrow0$ and $u^k\rightarrow u$ such that for any $y^k\in S(\bar x+t_ku^k)$, $y^k\notin k\mathbb B$. Obviously for both of the two cases, $\{y^k\}$ does not contain a convergent subsequence, which contradicts the directional inner semi-compactness. Hence the directional inner semi-compactness implies the directional inf-compactness. 
Note that our definition of directional inner semi-compactness condition differs from the one
  proposed by \cite{Long}. In fact the one defined in \cite{Long} can not recover its non-directional counterpart when the direction is equal to the origin; see e.g., \cite[Definition 1.63]{Aub2} for the definition of inner semi-compactness. The inner semi-compactness defined in Definition \ref{ic} is weak and, when considering the direction 0, coincides with the classical one.
}

The following definition gives a directional version of the classical inner semi-continuity (see e.g. \cite[Definition 1.63]{Aub2}).
\begin{defn}[Directional Inner Semi-continuity] Given $\bar y\in S(\bar x)$, we say that  the optimal solution map $S(x)$ is inner semi-continuous at $(\bar x,\bar y)$ in direction $u$, if for any sequences $ t_k\downarrow 0, u^k\rightarrow u$, there exists a sequence $y^k\in S(\bar x+t_ku^k)$ converging to $\bar 
	y$. When $u=0$ in the above, we say that $S(x)$ is inner semi-continuous at $(\bar x,\bar y)$.
\end{defn}
Note that \cite[Definition 4.4(i)]{Long}  introduced a directional inner semi-continuity  which requires $y^k\xrightarrow{v} \bar y$ for some $v$. Since $y^k\xrightarrow{v} \bar y$ implies that  $y^k \rightarrow \bar y$, their directional inner semicontinuity  is stronger than ours.
Given a direction we define a subset of the solution $S(\bar x)$ as below. It coincides with the solution set when $u=0$ and may be strictly contained in the solution set if the direction $u$ is nonzero. {Furthermore, one can easily obtain that the directional inner semi-continuity is stronger than the directional inner semi-compactness.}
\begin{defn}[Directional Solution]
	The optimal solution in direction $u$ is the set defined by
	$$S(\bar x;u)=\{y\in S(\bar x)| \exists t_k\downarrow 0, u^k\rightarrow u, y^k\rightarrow y, y^k\in S(\bar x+t_k u^k)\}.$$  
\end{defn}
It is obvious that if the optimal solution map $S(x)$ is inner semi-continuous at $(\bar x,\bar y)\in gph S$ in direction $u$, then $\bar y\in S(\bar x; u)$. Also if $(P_x)$ is upper stable at $\bar x$ in direction $u$ in the sense of Janin (see \cite[Definition 3.4]{Janin}), then $S(\bar x;u)\neq\emptyset$.

\begin{defn}[RCR regularity](\cite[Definition 1]{MS}) 
We say that {the the feasible  map ${\cal F}(x)$  is relaxed constant rank (RCR) regular at} $(\bar x,\bar y ) \in {\rm gph }{\cal F}$
 if there exists $\delta>0$ such that for any index subset $K\subseteq I_g(\bar x,\bar y)$, the family of gradient vectors $\nabla_y g_j(x,y), j\in K$, has the same rank at all points $(x,y) \in \mathbb{B}_\delta(\bar x,\bar y)$.
\end{defn}

We now define a directional version of the Robinson Stability \cite[Definition 1.1]{HM}).
\begin{defn}[Directional Robinson stability]
We say that the feasible map ${\cal F}(x)$ satisfies Robinson stability ${\rm (RS)}$ property at $(\bar x,\bar y)\in {\rm gph }{\cal F}$ in direction $u\in\mathbb R^n$ if there exist positive scalars $\kappa, \epsilon, \delta$ such that 
\begin{equation}\label{RS}
{\rm dist}(y, {\cal F}(x))\leq \kappa \|g_+(x,y)\| \quad \forall x\in \bar x+{\cal V}_{\epsilon,\delta}(u), y\in  \mathbb{B}_\epsilon(\bar y).
\end{equation}
\end{defn}
If RS holds at $(\bar x,\bar y)$ in direction $u=0$, we say that  RS holds at $(\bar x,\bar y)$ (\cite[Definition 1.1]{HM}). 
RS means that the error bound condition holds at $\bar y$ uniformly in a directional neighborhood of $\bar x$ in direction $u$.
Note that  RS in direction $u$ is equivalent to R-regularity with respect to set $\bar x+{\cal V}_{\epsilon,\delta}(u)$ as defined in \cite[Definition 2]{MS}.

\begin{prop}[Sufficient Conditions for RS] 
	If the system $g(x,y)\leq0$ satisfies one of the following conditions at $(\bar x,\bar y)$, then RS holds at $(\bar x,\bar y)$.
	\begin{itemize}
		\item[1.] {$g(x,y)=Ax+By+c$, where $A,B$ are $n\times p$ and $m\times p$ matrices, respectively,
	and the feasible region $\mathcal F(x)$ is nonempty near $\bar x$.}
		\item[2.]  The set $\{\nabla_y g_i(\bar x,\bar y)| i\in I_g(\bar x,\bar y)\}$ is linearly independent.
		\item[3.]   There exists no nonzero vector $\lambda\in\mathbb R^p_+$ such that $\lambda\perp g(\bar x,\bar y)$ and $\nabla_y g(\bar x,\bar y)^T\lambda=0$.
	\end{itemize}
\end{prop}
\beginproof By condition 1, for any $x$ close to $\bar x$, the feasible region of the linear inequality system $\mathcal F(x)$ is nonempty.
Then by Hoffman's lemma \cite{Hoffman}, the error bound holds, i.e., there exists $\kappa>0$ such that
 \begin{equation}\label{Hoffman}
 {\rm dist}(y,\mathcal F(x))\leq\kappa \|(Ax+By+c)_+\|\quad \mbox{for }\forall y.
 \end{equation}  Moreover since by \cite{Hoffman}, the $\kappa$ depends only on matrix $B$,  it is independent of variable $x$ and inequality (\ref{Hoffman}) holds at any $(x,y)$ near $(\bar x,\bar y)$. Hence, condition 1 implies RS. Condition 2 implies condition 3 and it is well known that condition 3  is equivalent to saying that  the  condition ${\rm dist}(y, {\cal F}(x))\leq \kappa \|g_+(x,y)\| \quad \forall (x,y)\in  \mathbb{B}_\epsilon(\bar x, \bar y)$ holds  for sufficiently small $\epsilon>0$.
\endproof
One can refer to \cite{MS, GLYZ, HM} and the references therein for more sufficient conditions for RS. {The following proposition shows that   the directional RS implies the directional MSCQ.

\begin{prop}\label{relation}
Suppose that the feasible map $\mathcal F$ satisfies RS at $(\bar x,\bar y)\in {\rm gph }{\cal F}$ in direction $u$. Then the system $g(x,y)\leq0$ satisfies the directional MSCQ at $(\bar x,\bar y)$ in direction $(u,v)$ for any $v\in\mathbb R^m$.
\end{prop}
\beginproof
Since RS for $\mathcal F$ holds at $(\bar x,\bar y)$ in direction $u$, i.e., there exist numbers $\kappa>0,\epsilon>0,\delta>0$ such that
\[
{\rm dist}(y,\mathcal F(x))\leq\kappa \|g_+(x,y)\|,
\]  
for all $x\in\bar x+{\cal V}_{\epsilon,\delta}(u)$ and $y\in \mathbb B_\epsilon(\bar y)$.
Then we obtain
\[
{\rm dist}((x,y),g^{-1}(\mathbb R^p_-))={\rm dist}((x,y),gph\mathcal F)\leq dist(y,\mathcal F(x))\leq\kappa \|g_+(x,y)\|,
\]  
for all $x\in\bar x+{\cal V}_{\epsilon,\delta}(u)$ and $y\in \mathbb B_\epsilon(\bar y)$. This means that the directional MSCQ holds for the system $g(x,y)\leq0$ at $(\bar x,\bar y)$ in direction $(u,v)$ for any $v\in\mathbb R^m$.
\endproof

Recall that the lower Dini directional derivative of the feasible  map $\mathcal F(x)$ at a point $(\bar x,\bar y)\in {\rm gph}\mathcal F$ in direction $u$ is defined as
\[
D_+\mathcal F(\bar x,\bar y;u):=\liminf_{t\downarrow0}\frac{\mathcal F(\bar x+tu)-\bar y}{t}=\{v| \exists o(t) \mbox{ s.t. } \bar y +tv +o(t)\in {\cal F} (\bar x+ tu)\}.
\]
Define the $y$-projection of the linearization cone of gph${\cal F}$ at $(\bar x,\bar y)$ in direction $u$, i.e.,
$$\mathbb L(\bar x,\bar y;u):=\{v\in\mathbb R^m|\nabla g_i(\bar x,\bar y)(u,v)\leq0, i\in I_g(\bar x,\bar y)\}.$$ By definition,  one always has $D_+\mathcal F(\bar x,\bar y;u) \subseteq \mathbb{L}(\bar x,\bar y;u)$.  
Since the directional MPEC R-regularity introduced in \cite[Lemma 3.3]{GLYZ} is weaker than our directional RS and $(P_x)$ is a special case of the problem studied in \cite{GLYZ} when the equilibrium constraints are omitted, the following results follow from \cite[Lemmas 3.3, 3.5]{GLYZ}.
\begin{lemma}\cite[Lemmas 3.3, 3.5]{GLYZ}\label{equal}
	Let $\bar y\in \mathcal F(\bar x)$. Suppose
	 either the feasible map $\mathcal F$ satisfies RS at $(\bar x,\bar y)$ in direction $u$
	 or $D_+\mathcal F(\bar x,\bar y;u)\neq\emptyset$ and $\mathcal F$ is RCR-regular at $(\bar x,\bar y)$.
 Then $D_+\mathcal F(\bar x,\bar y;u) = \mathbb{L}(\bar x,\bar y;u).
$
\end{lemma}

{The following results will be needed in Proposition \ref{ddv} and Theorem \ref{Lip}.} 
\begin{lemma}\label{dcont}Suppose that the restricted inf-compactness holds at $\bar x$ in direction $u$  with compact set $\Omega_u$  and there exists $\bar y\in S(\bar x)$ such that ${\cal F}$ satisfies RS at $(\bar x,\bar y)$ in direction $u$. Then $D_+\mathcal F(\bar x,\bar y;u)=\mathbb L(\bar x,\bar y;u)\neq\emptyset$, $\displaystyle \bar y\in\liminf_{x\xrightarrow{u}\bar x}\mathcal F(x)$  and $S(\bar x;u)\neq\emptyset$.
	And for any $l>0, \exists\epsilon, \delta>0$ such that for $\forall x\in\bar x+ V_{\epsilon,\delta}(u),\ \exists y\in S(x)\cap\Omega_u$ satisfying  ${\rm dist}(y,S(\bar x;u)\cap\Omega_u)<l$.
\end{lemma}
\beginproof  {Since RS holds at $(\bar x,\bar y)$ in direction $u$, by Lemma \ref{equal}, $D_+\mathcal F(\bar x,\bar y;u)=\mathbb L(\bar x,\bar y;u)$.} 
Moreover by $(\ref{RS})$ there exist positive scalars $\kappa, \epsilon, \delta$, such that for any $x\in\bar x+{\cal V}_{\epsilon,\delta}(u)$,  \begin{equation}\label{isocalm}
{\rm dist}(\bar y,\mathcal F(x))\leq\kappa \|g_+(x,\bar y)\|\leq\kappa\|g(x,\bar y)-g(\bar x,\bar y)\|\leq L_g\kappa\|x-\bar x\|,
\end{equation}
where $L_g>0$ is the Lipschitz modulus of $g(x,\bar y)$ around $\bar x$. Then for any sequences $t_k\downarrow0, u^k\rightarrow u$, by $(\ref{isocalm})$ we can find a sequence $y^k\in\mathcal F(\bar x+t_ku^k)$ such that $\|\bar y-y^k\|\leq L_g\kappa\|\bar x+t_ku^k-\bar x\|$, {which implies that} $y^k\rightarrow\bar y$. By Definition \ref{innerlimit}, this means that $\displaystyle{\bar y\in\liminf_{x\xrightarrow{u}\bar x}\mathcal F(x)}$. 
Since  $\{\frac{y^k-\bar y}{t_k}\}$ is bounded, taking a subsequence if necessary, we can find $v\in\mathbb R^m$ such that $v^k:=\frac{y^k-\bar y}{t_k}$ converges to $v$. {Since for each $i\in I_g(\bar x,\bar y)$, $g_i(\bar x+t_ku^k,y^k)\leq0$,} it follows that  $v\in\mathbb L(\bar x,\bar y;u)$.
{We also have $\limsup_k V(\bar x+t_ku^k) \leq \lim_k f(\bar x+t_ku^k,y^k)=V(\bar x)$. It follows that
since the restricted inf-compactness holds at $\bar x$ in direction $u$,  for each $k$ large enough, there exists $\tilde y^k\in S(\bar x+t_ku^k)\cap\Omega_u$.  By the compactness of $\Omega_u$, the sequence $\{\tilde y^k\}$ is bounded. Without loss of generality, assume $\tilde y:=\lim_k\tilde y^k\in\Omega_u$.
Since
\begin{eqnarray*}
&& f(\bar x,\tilde y)=\lim_{k}f(\bar x+t_ku^k,\tilde y^k)=\lim_kV(\bar x+t_ku^k)\leq\lim_{k}f(\bar x+t_ku^k,\bar y+t_kv^k)=f(\bar x,\bar y)=V(\bar x),\\
&& g(\bar x, \tilde y) = \lim_{k} g(\bar x+t_ku^k, \tilde y^k) \leq 0
\end{eqnarray*}
 we can obtain $\tilde y\in S(\bar x)$. Consequently, $\tilde y \in S(\bar x;u)\cap\Omega_u$.

We prove the last statement by contradiction. Assume there exist $l>0$ and for $k$ large enough, $\bar x+t_ku^k\in \bar x+{\cal V}_{\frac{1}{k},\frac{1}{k}}(u)$ and $\tilde y^k\in S(\bar x+t_ku^k)\cap\Omega_u$ such that  ${\rm dist}(\tilde y^k,S(\bar x;u)\cap\Omega_u)\geq l$. Taking the limit as $k\rightarrow\infty$, ${\rm dist}(\tilde y,S(\bar x;u)\cap\Omega_u)\geq l$, which contradicts $\tilde y\in S(\bar x;u)\cap\Omega_u$. The proof is complete. 
\endproof

\subsection{Directional derivative of the value function}

In this subsection, we study the directional differentiability of the value function. 
{In} the following proposition we derive the formula for the directional derivative of the value function. Our result improves the correspongding classical result in \cite[Theorem 3.11]{GLYZ} in that in the formula the directional solution instead of the solution set is used and
the NNAMCQ holding at each $y\in S(\bar x)$ is replaced by the directional RS which is in general weaker.

\begin{prop}\label{ddv}
Assume that ${\cal F}$ is RCR-regular at each  $(\bar x, y)\in{\rm gph} S$. Suppose that the restricted inf-compactness holds at $\bar x$ in direction $u$ and RS is satisfied at each $(\bar x, y)\in{\rm gph} S$ in direction $u$.  Then the value function is directionally differentiable at $x=\bar x$ in direction $u$ and		
\begin{equation}
	V'(\bar x;u)=\min_{y\in  S(\bar x;u)}\min_{v\in \mathbb{L}(\bar x, y;u)}\nabla f(\bar x, y)(u,v)=\min_{y\in  S(\bar x;u)}\max_{\lambda\in\Lambda(\bar x,y)}\nabla_x{\cal L}(\bar x,y;\lambda)u.\label{equality}
\end{equation}
\end{prop}
\beginproof
Since ${\cal F}$ satisfies RS at each $(\bar x, y)\in{\rm gph} S$ in direction $u$ {and the restricted inf-compactness holds at $\bar x$ in direction $u$}, by Lemma \ref{dcont}, $S(\bar x;u)\neq\emptyset$ and $D_+{\cal F}(\bar x,y;u)\neq\emptyset$ for $\forall y\in S(\bar x;u)$.
Then for any given $y\in S(\bar x;u)$, there is $v\in D_+\mathcal F(\bar x, y;u)$. It follows that   there exists $o(t)$ such that $ y+tv+o(t)\in\mathcal F(\bar x+tu)$ for $t\geq0$. Thus 
we have
\begin{eqnarray}\label{udini}
V_+'(\bar x;u)
:=\limsup_{t\downarrow0}\frac{V(\bar x+tu)-V(\bar x)}{t} 
&\leq &\limsup_{t\downarrow0}\frac{f(\bar x+tu, y+tv+o(t))-f(\bar x, y)}{t}\nonumber \\
&=&\nabla f(\bar x, y)(u,v). \label{upestimate}
\end{eqnarray}
On the other hand, let $t_k\downarrow0$ be the sequence satisfying
\begin{equation*}
V_-'(\bar x;u):=\liminf_{t\downarrow0}\frac{V(\bar x+tu)-V(\bar x)}{t}=\lim_{k\rightarrow \infty}\frac{V(\bar x+t_ku)-V(\bar x)}{t_k}.
\end{equation*}
By $(\ref{udini})$, for any $\epsilon>0$ and any sequence $t_k\downarrow0$,  $V(\bar x+t_ku)< V(\bar x)+\epsilon$ for $k$ large enough. Since the restricted inf-compactness holds at $\bar x$ in direction $u$ with a compact set $\Omega_u$, there exists a sequence $y^k\in S(\bar x+t_ku)\cap\Omega_u$ for $k$ large enough. Without loss of generality, define $\tilde y:=\lim_ky^k$. Then 
\begin{eqnarray*}
&& f(\bar x,\tilde y)=\lim_{k\rightarrow \infty}f(\bar x+t_ku,y^k)=\lim_{k\rightarrow \infty}{ V}(\bar x+t_ku)\leq { V}(\bar x),\\
&& g(\bar x,\tilde y)= \lim_{k\rightarrow \infty} g(\bar x+t_ku,y^k) \leq 0.
\end{eqnarray*}  {This means $\tilde y\in S(\bar x)\cap\Omega_u$. Moreover it is clear that
$\tilde y\in S(\bar x;u)\cap\Omega_u.
$ }

Since ${\cal F}$ is RCR regular at each $y\in S(\bar x;u)$ and  $ D_+\mathcal F(\bar x, y;u)\not =\emptyset$, {by Lemma \ref{equal}} we have
\begin{equation}D_+{\cal F}(\bar x, y; u)=\mathbb{L}(\bar x, y; u) \quad \forall y\in S(\bar x;u).\label{lowerd} \end{equation}
Moreover by \cite[Lemma 5]{MS}, for sufficiently large $k$, there exist $\kappa>0$ independent of $k$ and a sequence $\bar y^k\in\mathcal F(\bar x)$ such that
\begin{equation*}
\|y^k-\bar y^k\|\leq\kappa\|\bar x+t_ku-\bar x\|,\ g_j(\bar x+t_ku,y^k)-g_j(\bar x,\bar y^k)\leq0,\ j\in  I_g(\bar x,\tilde y).
\end{equation*}
Consequently, $\{\frac{y^k-\bar y^k}{t_k}\}$ is bounded. Taking a subsequence if necessary, we assume that  $\tilde v:=\lim_{k\rightarrow \infty} \frac{y^k-\bar y^k}{t_k}$ and then $y^k=\bar y^k+t_k \tilde v+o(t_k)$.
Thus, we obtain $\nabla g_i(\bar x,\tilde y)(u,\tilde v)\leq0, i\in  I_g(\bar x,\tilde y)$. This implies that 
$
\tilde v\in \mathbb{L}(\bar x,\tilde y;u).
$
Furthermore, since $\bar y^k\in\mathcal F(\bar x)$, we have
\begin{eqnarray}
\nonumber
V'_-(\bar x;u)&&=\lim_{k\rightarrow \infty}\frac{V(\bar x+t_ku)-V(\bar x)}{t_k}\\
\nonumber
&&\geq\lim_{k\rightarrow\infty}\frac{f(\bar x+t_ku,y^k)-f(\bar x,\bar y^k)}{t_k}\\
\nonumber
&&=\lim_{k\rightarrow\infty}\frac{f(\bar x+t_ku,\bar y^k+t_k \tilde v+o(t_k))-f(\bar x,\bar y^k)}{t_k}\\
&&=\nabla f(\bar x,\tilde y)(u,\tilde v).\label{low}
\end{eqnarray}
It follows that 
\begin{equation}\label{loestimate}
V_-'(\bar x;u)\geq\nabla f(\bar x, \tilde y)(u, \tilde v)\geq {\min_{ y\in S(\bar x;u)\cap\Omega_u}}\inf_{v\in  \mathbb{L}(\bar x, y;u)}\nabla f(\bar x, y)(u,v).
\end{equation}
Since  (\ref{upestimate}) holds for any $ y \in S(\bar x;u)\subseteq S(\bar x)$ and $v\in D_+\mathcal F(\bar x, y;u)=\mathbb{L}(\bar x, y; u)$,  where the equality follows from (\ref{lowerd}), we have
{\begin{equation} V_+'(\bar x;u)\leq \inf_{ y\in S(\bar x;u)}\inf_{v\in  \mathbb{L}(\bar x, y;u)}\nabla f(\bar x, y)(u,v)\leq\min_{ y\in S(\bar x;u)\cap\Omega_u}\inf_{v\in  \mathbb{L}(\bar x, y;u)}\nabla f(\bar x,                                                          y)(u,v) .\label{loestimatenew}
\end{equation}}
$(\ref{loestimate})$ and $(\ref{loestimatenew})$ imply that 
$$V_-'(\bar x;u)\geq\inf_{ y\in S(\bar x;u)}\inf_{v\in  \mathbb{L}(\bar x, y;u)}\nabla f(\bar x,                                                          y)(u,v)\geq V_+'(\bar x;u).$$ Hence 
$V(x)$ is directionally differentiable  at $\bar x$ in direction $u$ with the first equality in (\ref{equality}) holds.
{And the minimum with respect to $y$ in (\ref{equality}) can be attained on the set $S(\bar x;u)\cap\Omega_u$.
By the linear programming duality theorem, the second equality in (\ref{equality}) holds and the minimum with respect to $v$ can be attained.}
\endproof


In general, according to Proposition \ref{ddv}, one needs to ensure both RS and RCR regularity for the existence of the directional derivative. However, the following proposition shows that, if the solution set $S(x)$ is inner semi-continuous, only RCR-regularity is needed.
\begin{prop}\label{Semic} Suppose that the solution set $S( x)$ is inner semi-continuous at $(\bar x,\bar y)\in {\rm gph} S$ in direction $u$. Moreover assume that $\mathcal F$ is  RCR-regular at $(\bar x,\bar y)$. Then the value function $V(x)$ is  directionally differentiable at $\bar x$ in direction $u$ and
	$$
	V'(\bar x;u)=\min_{v\in \mathbb{L}(\bar x, \bar y;u)}\nabla f(\bar x, \bar y)(u,v)=\max_{\lambda\in\Lambda(\bar x,\bar y)}\nabla_x{\cal L}(\bar x,\bar y;\lambda)u.$$	
\end{prop}
\beginproof
Since $S( x)$ is inner semi-continuous at $(\bar x,\bar y)\in {\rm gph} S$ in direction $u$, we have that  the restricted inf-compactness holds at  $\bar x$ in direction $u$.

Next, we claim that since the directional inner semi-continuity and RCR-regularity hold, RS holds at $(\bar x,\bar y)$ in direction $u$. 
We prove this by contradiction. Assume RS does not hold at $(\bar x,\bar y)$ in direction $u$. Then there exist sequences $x^k\xrightarrow{u} \bar x$ and $y^k\rightarrow \bar y$ satisfying that
\begin{equation}\label{contradiction}
{\rm dist}(y^k,\mathcal F(x^k))> k\|g_+(x^k,y^k)\|.
\end{equation}

{Since $S(x)$ is inner semi-continuous at $(\bar x,\bar y)$ in direction $u$, we have
	\[
	\bar y\in\liminf_{x\xrightarrow{u}\bar x} S(x)\subseteq\liminf_{x\xrightarrow{u}\bar x}\mathcal F(x).
	\]}
Then for sufficiently large $k$ there exists a sequence $\tilde y^k\in \mathcal F(x^k)$ such that $\tilde y^k\rightarrow\bar y$. Let $\bar y^k$ be the projection of $y^k$ onto $\mathcal F(x^k)$. We obtain
\[
\|y^k-\bar y^k\|\leq\|y^k-\tilde y^k\|\rightarrow 0\ \mbox{as}\ k\rightarrow\infty.
\]
Then following the proof of \cite[Lemma  3.5]{GLYZ} for the case when the number of complementarity constraints is 0, we can find some scalar $M>0$ such that
\begin{equation*}
{\rm dist}(y^k,\mathcal F(x^k))\leq M\|g_+(x^k,y^k)\|
\end{equation*} contradicting (\ref{contradiction}).
Hence, the assumption is false and RS for $\mathcal F$ holds at $(\bar x,\bar y)$ in direction $u$. 

Now RCR holds at $(\bar x,\bar y)$ and RS holds at $(\bar x,\bar y)$ in direction $u$.  By definition of the directional  inner semi-continuity of $S(x)$, we can always choose $\tilde y=\bar y$ in the proof of Proposition \ref{ddv}. Hence the result follows from Proposition \ref{ddv}.
\endproof
\subsection{Directional Lipschitz continuity of the value function}
In this subsection we study  sufficient conditions for the directional Lipschitz continuity of $V(x)$. 

The classical {criterion} for  guaranteeing the Lipschitz continuity of the value function, is a combination of the uniform compactness condition and  {MFCQ} holding at each $y\in S(\bar x)$, see e.g. \cite[Theorem 5.1]{Gauvin}.
	The following theorem gives sufficient conditions for the directional Lipschitz continuity of the value function under weaker assumptions.} {When $u=0$, it recovers the result in \cite[Theorem 5.5]{BMR}.}


\begin{thm}\label{Lip}
\begin{itemize}
\item [(i)]
	Suppose that the restricted inf-compactness holds at $\bar x$ in direction $u$ with compact set $\Omega_u$ and the feasible map ${\cal F}(x)$ satisfies RS at each $(\bar x, y)\in {\rm gph} S$.
	Then $V(x)$ is directionally Lipschitz continuous at $\bar x$ in direction $u$. 
\item[(ii)] 
	Suppose there exists $\bar y\in S(\bar x)$ such that $S(x)$ is inner semi-continuous at $(\bar x,\bar y)$ in direction $u$, and the feasible map ${\cal F}(x)$  satisfies RS at $(\bar x, \bar y)$ in direction $u$.
	Then $V(x)$ is directionally Lipschitz continuous at $\bar x$ in direction $u$. 
\end{itemize}	Furthermore, if $(i)$ or $(ii)$ holds in direction $u=0$ then $V(x)$ is Lipschitz at $\bar x$.
\end{thm}
\beginproof 
Since $\mathcal F$ satisfies RS at each $(\bar x,y)\in{\rm gph} S$ in direction $u$, by Lemma \ref{dcont}, $S(\bar x;u)\cap\Omega_u\neq\emptyset$. And by the compactness of $\Omega_u$ and Borel-Lebesgue covering theorem, there exist positive scalars $\epsilon,\delta,\kappa$ such that 
\begin{equation}\label{Borel-Lebesque}
{\rm dist}(y,\mathcal F(x))\leq\kappa \|g_+(x,y)\| \quad \forall x\in \bar x+{\cal V}_{\epsilon,\delta}(u), y\in (S(\bar x;u)\cap\Omega_u)+\epsilon\mathbb B.
\end{equation}

By Lemma \ref{dcont}, choosing $\epsilon, \delta$ small enough, we have for any $x,x'\in\bar x+{\cal V}_{\epsilon,\delta}(u)$, there exist $y\in S(x)\cap\Omega_u,\ y'\in S(x')\cap\Omega_u$ close enough to $S(\bar x;u)\cap \Omega_u$. 
Without loss of generality  assume
$x,x' \in \bar x+{\cal V}_{\epsilon,\delta}(u)$ and $y,y' \in (S(\bar x;u)\cap\Omega_u)+\epsilon\mathbb B$.  Then by (\ref{Borel-Lebesque}) we can find $\bar y\in \mathcal F(x),\ \bar y'\in \mathcal F(x')$ such that
\begin{align*}
&\|y-\bar y'\|\leq\kappa\|g(x',y)-g(x,y)\|\leq2\kappa\|\nabla_xg(x,y)\|\|x-x'\|,\\
&\|y'-\bar y\|\leq\kappa\|g(x,y')-g(x',y')\|\leq2\kappa\|\nabla_xg(x',y')\|\|x-x'\|.
\end{align*}
Since $\nabla_xg(x,y)$ is continuous and $\{\bar x+{\cal V}_{\epsilon,\delta}(u)\}\times \left(S(\bar x;u)\cap\Omega_u\right)$ is bounded, by Weirstrass extreme value theorem, there exists a positive scalar $M$ such that $2\kappa\|\nabla_xg(x,y)\|\leq M$ for any $(x,y)\in\{\bar x+{\cal V}_{\epsilon,\delta}(u)\}\times \left(S(\bar x;u)\cap\Omega_u\right)$. Similarly, since $\nabla f(x,y)$ is continuous, hence, locally bounded. Choosing $M'$ large enough, we have 
\begin{align*}
&\|f(x,y)-f(x',\bar y')\|\leq M'\|(x,y)-(x',\bar y')\|\leq M'(1+M)\|x-x'\|,\\
&\|f(x,\bar y)-f(x',y')\|\leq M'\|(x,\bar y)-(x',y')\|\leq M'(1+M)\|x-x'\|.
\end{align*} Then since
$f(x,y)-f(x',\bar y')\leq V(x)-V(x')\leq f(x,\bar y)-f(x',y')$, we have $$\|V(x)-V(x')\|\leq\max\{\|f(x,y)-f(x',\bar y')\|,\|f(x,\bar y)-f(x',y')\|\}\leq M'(1+M)\|x-x'\|.$$ This means $V(x)$ is directionally Lipschitz continuous at $\bar x$ in direction $u$ and (i) is proved.

Next, we prove (ii). If there exists $\bar y\in S(\bar x)$ such that $S(x)$ is inner semi-continuous at $(\bar x,\bar y)$ in direction $u$, $\bar y\in S(\bar x;u)\neq\emptyset$ and the restricted inf-compactness holds at $\bar x$ in direction $u$. Then one can easily replace $S(\bar x;u)\cap\Omega_u$ by $\{\bar y\}$ in the proof above and obtain the directional Lipschitz continuity of $V(x)$ under RS at $(\bar x,\bar y)$ in direction $u$.
\endproof
        \subsection{Directional subdiffentials of the value function}
In this subsection, we study the directional subdifferential of the value function of $(P_x)$. 
{First, 
	 we derive an upper estimate for the  directional subdifferential of the value function  in terms of the problem data.}  For any $x,y,u$, suppose $V'(x;u)$ exists. We denote by 
	\begin{equation}\Sigma(x, y, u):= \{v\in \mathbb L( x, y; u)| V'( x;u)=\nabla f(x, y)(u,v)\}.
	\label{Sigma}
\end{equation}
\begin{thm}
\label{estimates}
Let $u\in \mathbb{R}^n$.

(i) Suppose that the restricted inf-compactness holds at $\bar x$ in direction $u$ with compact set $\Omega_u$. Suppose that $\mathcal F$ is RCR-regular at each $(\bar x,y)\in{\rm gph} S$ and satisfies RS at each $(\bar x, y)\in{\rm gph} S$ in direction $u$. Then $V(x)$ is directionally Lipschitz continuous at $\bar x$ in direction $u$ and
\begin{equation} \emptyset \not = \partial V(\bar x;u) \subseteq \Theta(\bar x;u) \label{est}\end{equation}
where 
	\begin{align}
& \Theta(\bar x;u) := \nonumber \\
	&\bigcup_{\tilde y\in  S(\bar x;u)\cap \Omega_u}\big (\bigcup_{v\in\Sigma(\bar x, \tilde y, u)}\left \{\nabla_x f(\bar x, \tilde y)+\nabla_x g(\bar x, \tilde y)^T\lambda_g|\right.
		  \left.\begin{array}{l}
	\lambda_g \in \Lambda(\bar x,\tilde y)\cap \{\nabla g(\bar x, \tilde y)(u,v) \} ^\perp \end{array} \right \} 
	  \big ). 
	\label{theta}
 	\end{align}
	
(ii)  Suppose that there exists $\bar y\in S(\bar x)$ such that $S(x)$ is inner semi-continuous at $(\bar x,\bar y)$ in direction $u$ and {${\mathcal F}$ is RCR-regular at $(\bar x,\bar y)$}, then $V(x)$ is directionally Lipschitz at $\bar x$ in direction $u$ and
		\begin{align*}
	\emptyset \not =\partial V(\bar x;u)\subseteq
	&\bigcup_{v\in \Sigma(\bar x, \bar y, u)} \left \{\nabla_x f(\bar x, \bar y)+\nabla_x g(\bar x, \bar y)^T\lambda_g\bigg|\begin{array}{l}
	\lambda_g \in \Lambda(\bar x,\bar y)\cap \{\nabla g(\bar x, \bar y)(u,v) \} ^\perp \end{array} \right \}.
 	\end{align*}
\end{thm} 
\beginproof 
(i) By Proposition \ref{ddv}, $V(x)$ is directional differentiable  at $\bar x$ in direction $u$. It follows that  for any sequence $\epsilon_k\downarrow0$, there is a sequence $t_k\downarrow0$ such that for $k$ large enough, we have $V(\bar x+t_ku)<V(\bar x)+\epsilon_k$. By the assumption of the directional restricted inf-compactness, for $k$ large enough, there exists $\hat y^k\in S(\bar x+t_ku)\cap\Omega_u$. This means that $\{\hat y^k\}$ is bounded. Without loss of generality, there exists $\hat y=\lim_k\hat y^k$. And we know that $f(\bar x,\hat y)=\lim_kV(\bar x+t_ku)\leq V(\bar x)$. Hence, $\hat y\in S(\bar x;u)\neq\emptyset$. 

{Since the directional restricted inf-compactness holds at $\bar x$ in direction $u$ and RS is satisfied at $(\bar x,y)$ in direction $u$ for each $y\in S(\bar x)$, by Theorem \ref{Lip}(i), $V(x)$ is directionally Lipschitz continuous at $\bar x$ in direction $u$. Then by the well-known Rademacher's Theorem and \cite[Theorem 9.13]{RW}, $\partial V(\bar x;u)\neq\emptyset$. 

Let $\zeta\in\partial V(\bar x;u)$. Then by definition, there exist sequences $t_k\downarrow0,\ u^k\rightarrow u,\ \zeta^k\rightarrow\zeta$ such that $V(\bar x+t_ku^k)\rightarrow V(\bar x)$
and 
$\zeta^k\in \widehat \partial V(\bar x+t_ku^k)$. It follows that 
$V(\bar x+t_ku^k)<V(\bar x)+\epsilon$ for all $k$ large enough and hence by  the directional restricted inf-compactness, there exists $y^k\in S(\bar x+t_ku^k)\cap \Omega_u$.  Passing to a subsequence if necessary, we may assume that $y^k\rightarrow\tilde y$. Hence, by the continuity of $f(x,y)$, $\tilde y\in S(\bar x;u)\cap \Omega_u$.

For each $k$, since $\zeta^k\in \widehat \partial V(\bar x+t_ku^k)$, there exists a neighborhood ${\cal U}^k$ of $\bar x+t_ku^k$ satisfying
\begin{equation*}
V(x)-V(\bar x+t_ku^k)-\langle\zeta^k,x-(\bar x+t_ku^k)\rangle+\frac{1}{k}\|x-(\bar x+t_ku^k)\|\geq0\ \forall x\in{\cal U}^k.
\end{equation*}
It follows from the fact $V(x)=\displaystyle \inf_{y} \{f(x,y)+\delta_{\mathbb R^p_-}(g(x,y))\}$
and $y^k\in S(\bar x+t_ku^k)$, that
\begin{equation*}
f(x,y)-\langle\zeta^k,x-(\bar x+t_ku^k)\rangle+\frac{1}{k}\|x-(\bar x+t_ku^k)\|+\delta_{\mathbb R^p_-}(g(x,y))\geq f(\bar x+t_ku^k,y^k),
\end{equation*}
for any $(x,y)\in{\cal U}^k\times\mathbb R^m$.
{Hence the function 
\begin{equation*}
\phi_k(x,y):=
f(x,y)-\langle\zeta^k,x-(\bar x+t_ku^k)\rangle+\frac{1}{k}\|x-(\bar x+t_ku^k)\|+\delta_{\mathbb R^p_-}(g(x,y))\end{equation*} attains its local minimum at $(x,y)=(\bar x+t_ku^k,y^k)$.
Thus, by the well known Fermat's rule and the sum rule (\cite[Exercise 10.10]{RW}),} 
\begin{equation}
0\in\nabla f(\bar x+t_ku^k,y^k)-(\zeta^k,0)+\frac{1}{k}\bar{\mathbb B}\times \{0\}+\partial (\delta_{\mathbb R^p_-}\circ g)(\bar x+t_ku^k,y^k).\label{bounded}
\end{equation}

{Since $\mathcal F$ satisfies RS at $(\bar x,\tilde y)$ in direction $u$ and $y^k\rightarrow\tilde y$,}  by Proposition \ref{relation},  the directional MSCQ holds for the system $g(x,y)\leq0$ at each $(\bar x+t_ku^k,y^k)$ for $k$ sufficiently large. Hence by \cite[Theorem 5]{HO2005} we have
\begin{equation*}
\partial (\delta_{\mathbb R^p_-}\circ g)(\bar x+t_ku^k,y^k)=N_{g^{-1}(\mathbb R^p_-)}(\bar x+t_ku^k,y^k)
\subseteq \nabla g(\bar x+t_ku^k,y^k)^TN_{\mathbb R^p_-}(g(\bar x+t_ku^k,y^k)).
\end{equation*}
Hence
\begin{equation}\label{boundednew}
0\in\nabla f(\bar x+t_ku^k,y^k)-(\zeta^k,0)+\frac{1}{k} \bar{\mathbb{B}}\times \{0\}+\nabla g(\bar x+t_ku^k,y^k)^TN_{\mathbb R^p_-}(g(\bar x+t_ku^k,y^k)).
\end{equation}} 
{
Since RCR-regularity holds at $(\bar x,\bar y)$ and $y^k\in S(\bar x+t_ku^k)$, by \cite[Lemma 5]{MS}, for sufficiently large $k$, there exist $\kappa>0$ independent of $k$ and a sequence $\bar y^k\in\mathcal F(\bar x)$ such that
\begin{equation}\label{property}
\|y^k-\bar y^k\|\leq\kappa\|\bar x+t_ku^k-\bar x\|,\ g_j(\bar x+t_ku^k,y^k)\leq g_j(\bar x,\bar y^k),\ j\in  I_g(\bar x,\bar  y).
\end{equation}
Then $I_g(\bar x+t_ku^k,y^k)\subseteq I_g(\bar x,\bar y^k)$ and
{by observing the formula of $N_{\mathbb R^p_-}(\cdot)$, one can easily get the following relationship} \begin{align*}
N_{\mathbb R^p_-}(g(\bar x+t_ku^k,y^k))&=N_{\mathbb R^p_-}(g(\bar x,\bar y^k))\cap[g(\bar x,\bar y^k)-g(\bar x+t_ku^k,y^k)]^\perp\\
&=N_{\mathbb R^p_-}(g(\bar x,\bar y^k))\cap[\frac{g(\bar x,\bar y^k)-g(\bar x+t_ku^k,y^k)}{t_k}]^\perp.
\end{align*}
Define $v^k:=\frac{y^k-\bar y^k}{t_k}$ and  $y^k=\bar y^k+t_kv^k$. By $(\ref{property})$, $\{v^k\}$ is bounded. Without loss of generality, there exists $v=\lim_kv^k$. Then $\lim_k\frac{g(\bar x+t_ku^k, y^k)-g(\bar x,\bar y^k)}{t_k}=\nabla g(\bar x,\tilde y)(u,v)$. By $(\ref{property})$, $\bar y^k\rightarrow\tilde y$ and $v\in\mathbb L(\bar x,\tilde y;u)$. Taking the limit in $(\ref{boundednew})$, we have
\begin{equation*}
0\in\nabla f(\bar x,\tilde y)-(\zeta,0)+\nabla g(\bar x,\tilde y)^T(N_{\mathbb R^p_-}(g(\bar x,\tilde y))\cap[\nabla g(\bar x,\tilde y)(u,v)]^\perp).
\end{equation*} We obtain the existence of $\lambda_g\in N_{\mathbb R^p_-}(g(\bar x,\tilde y))\cap[\nabla g(\bar x,\tilde y)(u,v)]^\perp$ such that $\zeta=\nabla_xf(\bar x,\tilde y)+\nabla_x g(\bar x,\tilde y)^T\lambda_g$.
Furthermore, since $y^k\in S(\bar x+t_ku^k)$ and $\bar y^k\in \mathcal F(\bar x)$, we have 
\[
\lim_k\frac{f(\bar x+t_ku^k,\bar y^k+t_kv^k)-f(\bar x,\bar y^k)}{t_k}\leq\frac{V(\bar x+t_ku^k)-V(\bar x)}{t_k}=V'(\bar x;u).
\]
Then by Proposition \ref{ddv}, we know $\nabla f(\bar x,\tilde y)(u,v)=V'(\bar x;u)$. Hence, $v\in\Sigma(\bar x,\tilde y;u)$. The proof is complete.}

(ii) When $S(x)$ is inner semi-continuous at some point $\bar y\in S(\bar x)$ in direction $u$, one can choose $\tilde y=\bar y$. {And from the proof of Proposition \ref{Semic}, since the inner semi-continuity and RCR hold, $\mathcal F(x)$ satisfies RS at $(\bar x,\bar y)$ in direction $u$. Then by Theorem \ref{Lip}(ii), $V(x)$ is directionally Lipschitz continuous at $\bar x$ in direction $u$. Consequently, the results follows similarly as the proof of (i).}
\endproof
\cite[Theoerems 5.10 and 5.11]{Long} also gave an upper estimate of the value function of constrained programs  in terms of the coderivatives of the constraint mapping $\cal F$ under a  stronger version of directional inner semicontinuity \cite[Definition 4.4(i)]{Long} of $S(x)$. {Our result cannot be obtained from \cite[Theoerems 5.10 and 5.11]{Long} and is in a more explicit form.}

The following theorem provides an estimate of the directional Clarke subdifferential of the value function which will be used in the necessary optimality condition for bilevel programs.
\begin{thm}\label{subdiff}
Under the assumptions of Theorem \ref{estimates}(i), we have 
	\begin{equation*}
	\partial^c V(\bar x;u)\subseteq co\bigcup_{\tilde y\in  S(\bar x;u)\cap \Omega_u} W(\bar x, \tilde y,u,v) \quad \forall v\in \Sigma(\bar x,\tilde y,u).
	\end{equation*}
Under the assumptions of Theorem \ref{estimates}(ii), we have $$\partial^c V(\bar x;u)\subseteq  W(\bar x, 
\bar y,u,v),\quad \forall v\in \Sigma(\bar x,\bar y,u).$$ Here \begin{eqnarray*}
		W( x, y, u,v)&:=& \left\{\nabla_x f(x, y)+\nabla_x g(x,y)^T\lambda_g\bigg|\begin{array}{l}
			\lambda_g\in \Lambda(x,y)\cap \{\nabla g( x, y)(u,v)\}^\perp
		\end{array}\right\}.
	\end{eqnarray*}
\end{thm}
\beginproof
Since $\partial^c V(\bar x;u)=co {\partial} V(\bar x;u)$, by Theorem \ref{estimates}, we only need to show that  
$W(x,y,u,v_1)=W(x,y,u,v_2)$
 for any $v_1,v_2 \in \Sigma ( x,  y,u)$.
Define 
{\[
K(x,y,u,v):=\{\lambda_g|\nabla_y f(x, y)+\nabla_y g(x, y)^T\lambda_g=0, 0\leq\lambda_g\perp g( x, y),
\lambda_g\perp\nabla g(x, y)(u,v)\}.
\]} It suffices to show that $K(x,y,u,v_1)=K(x,y,u,v_2)$
 for any $v_1,v_2 \in \Sigma ( x,  y,u)$.
Since {$\nabla_y f(x, y)+\nabla_y g(x,y)^T\lambda_g=0,$} we have $\lambda_g^T\nabla_yg( x, y)v_i=-\nabla_yf( x, y)v_i$ for $i=1,2$. And since $\nabla_yf( x, y)v_1=\nabla_yf( x, y)v_2=V'( x;u)-\nabla_x f( x, y)u$, $\lambda_g^T\nabla_yg( x, y)v_1=\lambda_g^T\nabla_yg( x, y)v_2$. Hence, $\lambda_g^T\nabla g( x, y)(u,v_1)=\lambda_g^T\nabla g( x, y)(u,v_2)$. This implies $K(x,y,u,v_1)=K(x,y,u,v_2)$.\endproof

\section{Necessary optimality conditions for  bilevel programs}

{The main purpose of this section is to apply Theorem \ref{dKKT} to problem (VP) and the result of the directional sensitivity analysis of the value functions in Section 4 to derive a sharp necessary optimality condition for (VP) under a weak and verifiable constraint qualification.}

Assume that the value function is directionally differentiable at $\bar x$ in direction $u$. Define the linearization cone of (VP) at $(\bar x,\bar y)$ by  
$$\mathbb{L}(\bar x,\bar y):=
\left \{(u,v)|\begin{array}{l}
\nabla f(\bar x,\bar y)(u,v)\leq  V'(\bar x;u), \\
\nabla g_i(\bar x,\bar y)(u,v)\leq 0\ \forall i\in I_{g}(\bar x,\bar y), {\nabla G_i}(\bar x,\bar y)(u,v)\leq 0\ \forall i\in I_{G}(\bar x,\bar y)\end{array} \right \}.
$$
Note that although $f(x,y)-V(x) \leq 0$ is an inequality, it is in fact an equality constraint by the definition of the value function. {Hence under the Abadie constraint qualification,  one always have $\nabla f(\bar x,\bar y)(u,v)\geq V'(\bar x;u) $ for all $(u,v)$ satisfying $\nabla g_i(\bar x,\bar y)(u,v)\leq 0\ i\in I_{g}(\bar x,\bar y)$.}  Therefore if the Abadie constraint qualification holds, in the linearization cone  the inequality  $\nabla f(\bar x,\bar y)(u,v)\leq V'(\bar x;u)$ can be equivalently replaced by the equality. 

Let $(\bar x,\bar y)$ be a feasible solution of (VP). Denote  the critical cone  of (VP) at $(\bar x,\bar y)$ by
\begin{align*}
C(\bar x,\bar y)
:=\{(u,v)\in \mathbb{L}(\bar x,\bar y)|&F(\bar x,\bar y)(u,v)\leq0\}.
\end{align*}

{We now apply Theorem \ref{dKKT} to (VP) and obtain the following necessary optimality   condition for  the bilevel program (BP).}
{\begin{thm}\label{opt0}
	Let $(\bar x,\bar y)$ be a local minimizer of {\rm (BP)}.
	Suppose that the value function $V(x)$ is directionally Lipschitz continuous and directionally differentiable at $\bar x$ in direction $u$ and  $(u,v)\in C(\bar x,\bar y)$.   Moreover suppose that (VP) is calm at $(\bar x,\bar y)$ in direction $(u,v)$.
	Then the directional KKT condition holds. That is, there exists  $({\lambda_V},\lambda_g,\lambda_G)$ such that
	\begin{eqnarray*}
		&& 0\in \nabla F(\bar x,\bar y)+\lambda_V \nabla f(\bar x,\bar y)+\lambda_V \partial (-V)(\bar x;u)\times \{0\} +
		\nabla g(\bar x,\bar y)^T\lambda_g+\nabla G(\bar x,\bar y)^T\lambda_G,\\
		&&\lambda_V\geq0,\ 0\leq\lambda_g\perp g(\bar x,\bar y),\ \lambda_g\perp\nabla g(\bar x,\bar y)(u,v),\ 0\leq\lambda_G\perp G(\bar x,\bar y),\ \lambda_G\perp\nabla G(\bar x,\bar y)(u,v).
	\end{eqnarray*}
\end{thm}}
Note that the form of the necessary optimality condition above must be different from the
one obtained by using the KKT approach or the combined approach since in our condition there are only first order derivatives of the lower level problem data involved while if the KKT condition was used then the second order derivatives of the  lower level problem data  must be involved. For the details of these kinds of comparisons, one is referred to \cite{HS}.

To apply Theorem \ref{opt0} to (VP), there are two issues to resolve. First, how to estimate $\partial(-V)(\bar x;u)$ in terms of the problem data of (VP)? Secondly, under what conditions does the directional Clarke calmness condition holds for (VP)?

Unfortunately, it is usually hard to estimate $\partial(-V)(x;u)$ directly and $\partial(-V)(x;u)\neq -\partial V(x;u)$. To this end, we use $\partial^c(-V)(x;u)$ as an upper estimate. By Proposition \ref{Prop2.4}, $\partial^c(-V)(x;u)=-\partial^cV(x;u)$ and Theorem \ref{subdiff} provides upper estimates for $\partial^cV(x;u)$.

Now we consider the second issue. The directional Clarke calmness condition, although weak, is implicit and so in general  is hard to verify. Naturally one hopes to find sufficient conditions for  the directional Clarke calmness which can be verified. Recall that for the discussion from Section \ref{Sec3}, for problem (VP),  the following implications hold:
\begin{center}
NNAMCQ\\
$\Downarrow$\\
FOSCMS in direction $(u,v)$\\
$\Downarrow$\\
Quasi-normality in direction $(u,v)$\\
$\Downarrow$\\
Metric subregularity in direction $(u,v)$\\
$\Downarrow$\\
Calmness in direction $(u,v)$.	
\end{center}


It is known \cite{YZ95} that classical constraint qualifications such as the NNAMCQ fails at each feasible point of (VP). A natual question is whether it is possible that the FOSCMS holds at a feasible point of (VP)? 

Next we try to answer this question. 
Given $(u,v)\in\mathbb R^{n+m}$,  define the single-valued map  $M_{(u,v)}:\ \mathbb R^n\times\mathbb R^m\rightarrow\mathbb{R}\times\mathbb R^p\times\mathbb R^q$ by
$$
M_{(u,v)}(x,y):=(f(x,y)-V(x)+\langle u,x-\bar{x}\rangle^3+\langle v,y-\bar{y}\rangle^3, g(x,y), G(x,y)).
$$

\begin{lemma}\label{Prop31} Let $(\bar x,\bar y)$ be a feasible point of {\rm (VP)}. Assume that the value function is directionally differentiable at $\bar x$ in any direction $u$ and 
	$ 0\not =(u,v) \in \mathbb{L}(\bar{x},\bar y)$.
	Then the system $M_{(u,v)}(x,y)\leq0$ does not satisfy the directional MSCQ at $(\bar{x},\bar y)$ in direction $(u,v)$.
\end{lemma}
\beginproof To concentrate on  the main idea we omit the upper level constraint $G(x,y) \leq 0$ in the proof. To the contrary, suppose that $M_{(u,v)}(x,y)\leq0$ satisfies the directional MSCQ at $(\bar{x},\bar y)$ in the nonzero direction $ (u,v) \in \mathbb{L}(\bar{x},\bar y)$. Then by definition of the directional MSCQ in direction $(u,v)$, $\exists \kappa>0$, for all sequences $t_k\downarrow0,\ u^k\rightarrow u,\ v^k\rightarrow v$, we have for sufficiently large $k$
\begin{align}\label{msr}
{\rm dist}&\left((\bar{x}+t_ku^k,\bar{y}+t_kv^k),\ {M^{-1}_{(u,v)}}(\mathbb R^{1+p+q}_-)\right)\notag\\
& \leq\kappa\ \|(M_{(u,v)}(\bar{x}+t_ku^k,\bar{y}+t_kv^k))_+\|\notag \\
&\leq \kappa \left((f(\bar{x}+t_ku^k,\bar{y}+t_kv^k)-V(\bar{x}+t_ku^k)+t_k^3\langle u,u^k\rangle^3+t_k^3\langle v,v^k\rangle^3)_+\right.\notag\\
&\left.+\|g_+(\bar x+t_ku^k,\bar{y}+t_kv^k)\|\right).
\end{align} 

Since $(u,v)\in L(\bar{x},\bar y)$, we have $g(\bar x,\bar y)+t_k\nabla g(\bar x,\bar y)(u,v)\leq0$. Hence,  
\begin{eqnarray*}
	&&\lim_{k\rightarrow\infty}\frac{\|g_+(\bar x+t_ku^k,\bar{y}+t_kv^k)\|}{t_k}\\
	&&\leq\lim_{k\rightarrow\infty}\frac{
		\|g(\bar{x}+t_ku^k,\bar{y}+t_kv^k)-t_k\nabla g(\bar x,\bar y)(u,v)-g(\bar{x},\bar y) \|}{t_k}=0.
\end{eqnarray*} 
Similarly, since $f(\bar x,\bar y)-V(\bar x)=0,$ we have $ \nabla f(\bar x,\bar y)(u,v)-V'(\bar x;u)\leq0$,
\begin{equation}\label{go}
\lim_{k\rightarrow \infty}\frac{(f(\bar{x}+t_ku^k,\bar{y}+t_kv^k)-V(\bar{x}+t_ku^k)+t_k^3\langle u,u^k\rangle^3+t_k^3\langle v,v^k\rangle^3)_+}{t_k}=0.
\end{equation}
Since for every $t_k>0$ sufficiently small, we can find a point $(x_{t_k},y_{t_k}) \in  {M^{-1}_{(u,v)}}(0,0)$  satisfying (\ref{msr}), then
\begin{equation}\label{contradict}
f(x_{t_k},y_{t_k})-V(x_{t_k})+\langle u, x_{t_k}-\bar{x}\rangle^3+\langle v, y_{t_k}-\bar{y}\rangle^3\leq 0.
\end{equation}
And by $(\ref{go})$,
$ \lim_{t_k\downarrow0}t_k^{-1}\|(\bar{x}+t_ku^k,\bar{y}+t_kv^k)-(x_{t_k},y_{t_k})\|=0$. 

Since $(u,v)\neq(0,0)$, by $(\ref{contradict})$ we have for every $k$ sufficiently large
\begin{eqnarray*}
	0&&\geq f(x_{t_k},y_{t_k})-V(x_{t_k})+\langle u, x_{t_k}-\bar{x}\rangle^3+\langle v, y_{t_k}-\bar{y}\rangle^3\\
	&&=
	f(x_{t_k},y_{t_k})-V(x_{t_k})+\langle u, \bar x+t_ku-\bar{x}\rangle^3+\langle v, \bar y+t_kv-\bar{y}\rangle^3+o(t_k^3)\\
	&&\geq f(x_{t_k},y_{t_k})-V(x_{t_k})+\frac{t_k^3}{2}(\|u\|^6+\|v\|^6)\\
	&&> f(x_{t_k},y_{t_k})-V(x_{t_k}),
\end{eqnarray*}
contradicting that $V(x_{t_k})=\displaystyle \inf_{y\in {\cal F}(x_{t_k})}f(x_{t_k},y)\leq f(x_{t_k},y_{t_k})$.
\endproof

We are now ready to give a negative answer on the question if the FOSCMS can be satisfied by a feasible solution of (VP). 
\begin{prop}\label{abnorm}
	Assume that the value function is directionally differentiable at $\bar x$ in any direction $u$ and $(u,v)\in C(\bar x,\bar y)$. Then  there exists  a nonzero vector $(\lambda,\mu,\nu)\in\mathbb R^{1+p+q}$ such that
	$\lambda\geq 0, 0\leq\mu\perp g(\bar x,\bar y),\ \mu\perp \nabla g(\bar x,\bar y)(u,v), 0\leq\nu\perp G(\bar x,\bar y),\ \nu\perp \nabla G(\bar x,\bar y)(u,v)$ and 
	\begin{equation}\label{vFOSCMS}
	0\in\lambda\partial(f- V)(\bar x,\bar y;(u,v))+\nabla g(\bar x,\bar y)^T\mu+\nabla G(\bar x,\bar y)^T\nu.
	\end{equation} 
	Hence FOSCMS fails at any feasible solution of (VP) in any critical direction.
\end{prop}
\beginproof 
Since by Lemma \ref{Prop31}, $M_{(u,v)}(x,y)\leq0$ does not satisfies the directional MSCQ at $(\bar x,\bar y)$ in direction $(u,v)$ and the directional MSCQ is weaker than FOSCMS, hence,
FOSCMS
must fail at $(\bar x,\bar y)$ in direction $(u,v)$.
By the sum rule in  \cite[Theorem 5.6]{Long} for the directional subdifferential,
\[
\partial(f(x,y)-V(x)+\langle u,x-\bar{x}\rangle^3+\langle v,y-\bar{y}\rangle^3)(\bar x,\bar y;(u,v))=\partial(f-V)(\bar{x},\bar y;(u,v)).
\] Hence by Definition \ref{qp}(2) the FOSCMS for the inequality system $M_{(u,v)}(x,y)\leq 0$ at $(\bar x,\bar y)$ is the same as the (\ref{vFOSCMS}) which means that FOSCMS for (VP) at $(\bar x,\bar y)$  in direction $(u,v)$ fails.
\endproof

{Since the directional quasi-normality is weaker than the FOSCMS, in the sequel, we try to apply the directional quasi-normality to (VP). Below, we combine Proposition \ref{qpms}, Lemma \ref{mscalm} and Theorem \ref{opt0} and obtain a sharp necessary optimality condition for (VP) under the directional quasi-normality.}
 
\begin{thm}\label{opt}
	Let $(\bar x,\bar y)$ be a local minimizer of {\rm (BP)}.
	Suppose that the value function $V(x)$ is directionally Lipschitz continuous and directionally differentiable at $\bar x$ in direction $u$ and  $(u,v)\in C(\bar x,\bar y)$.   Moreover suppose that   the directional quasi-normality holds at $(\bar x,\bar y)$ in direction $(u,v)$, i.e., there exists no nonzero vector $(\alpha,\nu_g,\nu_G)\in\mathbb R^{1+p+q}_+$ and \begin{eqnarray}
	&&	0\in \alpha \nabla f(\bar x,\bar y)-\alpha \partial^c V(\bar x;u)\times \{0\}+\nabla g(\bar x,\bar y)^T \nu_g
		+\nabla G(\bar x,\bar y)^T \nu_G,\label{dKKT1}\\
		&&\nu_g\perp g(\bar x,\bar y),\ \nu_g\perp\nabla g(\bar x,\bar y)(u,v),\  \nu_G\perp G(\bar x,\bar y),\ \nu_G \perp \nabla G(\bar x,\bar y)(u,v), \label{dKKT2}
	\end{eqnarray}
	and there exists sequences $t_k\downarrow0,\ (u^k,v^k)\rightarrow(u,v)$ such that 
	\begin{align}
		\alpha(f(\bar x+t_ku^k,\bar y+t_kv^k)-V(\bar x+t_ku^k))>0,\ &\mbox{if}\ \alpha>0,\label{sequencial1}\\
	g_i(\bar x+t_ku^k,\bar y+t_kv^k)>0,\ &\mbox{if}\ (\nu_g)_i> 0,  i\in I_g,  \nonumber  \\
		G_i(\bar x+t_ku^k,\bar y+t_kv^k)>0,\ &\mbox{if}\ (\nu_G)_i> 0, i\in I_G.\label{sequencial3}
	\end{align} Then the directional KKT condition holds. Moreover, there exists  $({\lambda_V},\lambda_g,\lambda_G)$ such that
	\begin{eqnarray*}
	&& 0\in \nabla F(\bar x,\bar y)+\lambda_V \nabla f(\bar x,\bar y)-\lambda_V \partial^c V(\bar x;u)\times \{0\} +
	\nabla g(\bar x,\bar y)^T\lambda_g+\nabla G(\bar x,\bar y)^T\lambda_G,\\
	&&\lambda_V\geq0,\ 0\leq\lambda_g\perp g(\bar x,\bar y),\ \lambda_g\perp\nabla g(\bar x,\bar y)(u,v),\ 0\leq\lambda_G\perp G(\bar x,\bar y),\ \lambda_G\perp\nabla G(\bar x,\bar y)(u,v).
	\end{eqnarray*}
\end{thm}
\beginproof
Define $\phi(x,y):=(f(x,y)-V(x),g(x,y),G(x,y))$ and $\lambda_\phi:=(\alpha,\nu_g,\nu_G)$. Then by assumption, $\phi (x,y)$ is directionally Lipschitz continuous and directionally differentiable at $(\bar x,\bar y)$ in direction $(u,v)$.  Since $(u,v)\in C(\bar x,\bar y)$, we have $\nabla\phi(\bar x,\bar y)(u,v)\leq0$. Then since $f(\bar x,\bar y)-V(\bar x)=0, g(\bar x,\bar y)\leq0, G(\bar x,\bar y)\leq0$, $(\ref{dKKT2})$ means $0\leq \lambda_\phi\perp \phi(\bar x,\bar y)$ and $\lambda_\phi\perp \nabla\phi(\bar x,\bar y)(u,v)$. Since  
\begin{eqnarray*}
\lefteqn{\partial (f-V)(\bar x,\bar y;(u,v))=\nabla f(\bar x,\bar y)+\partial(-V)(\bar x;u)\times\{0\} }   \\
&& \subseteq \nabla f(\bar x,\bar y)+\partial^c(-V)(\bar x;u)\times\{0\} \\
&& \subseteq \nabla f(\bar x,\bar y)-\partial^c V(\bar x;u)\times \{0\},
\end{eqnarray*} 
where the first equation follows from \cite[Theorem 5.6]{Long} and the second inclusion follows from Proposition \ref{Prop2.4}, (\ref{dKKT1})-(\ref{sequencial3}) imply that the directional quasi-normality defined in Definition \ref{qp} holds. Applying Proposition \ref{qpms}, Lemma \ref{mscalm} and Theorem \ref{opt0}, the proof is complete.  \endproof
{

	When the conditions in Proposition \ref{ddv} and Theorem  \ref{subdiff} hold, one can apply the formulas of $V'(\bar x;u)$ and the upper estimates for $\partial^cV(\bar x;u)$ obtained in section 4 and derive the directional KKT condition in terms of the problem data under the directional quasi-normality as below.
\begin{thm}
	Let $(\bar x,\bar y)$ be a local minimizer of {\rm (BP)} and $u\in \mathbb{R}^n$.
	Suppose that the feasible map ${\cal F}(x):=\{y|g(x,y)\leq 0\}$ is RCR-regular at each $(\bar x,y)\in gph S$ and   satisfies RS at each $(\bar x,y)\in gph S$ in direction $u$. Moreover assume that the restricted inf-compactness holds at $\bar x$ in direction $u$. Then the value function $V(x)$ is directionally Lipschitz continuous and directionally differentiable at $\bar x$ in direction $u$ with
	\begin{eqnarray*}
	V'(\bar x;u)&=&\min_{y\in  S(\bar x;u)}\max_{\lambda\in\Lambda(\bar x,y)}\nabla_x{\cal L}(\bar x,y;\lambda)u,\\
 \partial V(\bar x;u) &\subseteq &\Theta(\bar x;u),
	\end{eqnarray*}
	where $\Theta(\bar x;u)$ is defined as in (\ref{theta}). 
	Suppose that the directional quasi-normality holds at $(\bar x,\bar y)$ in direction $(u,v)\in C(\bar x,\bar y)$ in Theorem \ref{opt}, with $\partial^cV(\bar x;u)$ replaced by {$co(\Theta(\bar x;u))$}. Then the directional KKT condition {and} Theorem \ref{opt} holds with $\partial^cV(\bar x;u)$ replaced by {$co(\Theta(\bar x;u))$}.	
\end{thm}

When the the solution map $S(x)$ is directionally inner semi-continuous at the point of interest, the directional restricted inf-compactness holds and RCR-regularity implies RS. Consequently, we can obtain the directional quasi-normality condition and the KKT condition of (VP) in the following more verifiable forms.

\begin{thm}\label{Thm4.4}
	Let $(\bar x,\bar y)$ be a local minimizer of {\rm (BP)}.
  Suppose that the feasible map $\mathcal F(x)$ is RCR-regular at $(\bar x,\bar y)$ and $S(x)$ is inner semi-continuous at $(\bar x,\bar y)$ in direction $u$. {Then the value function is directionally Lipschitz continuous and directionally differentiable at $\bar x$ in direction $u$ and
 $V'(\bar x;u)=\max_{\lambda\in\Lambda(\bar x,\bar y)}\nabla_x{\cal L}(\bar x,\bar y;\lambda)u.$}
Suppose that there exists $v$ such that $(u,v)\in C(\bar x,\bar y)$.
Furthermore suppose that  there exists no nonzero vector $(\alpha,\nu_g,\nu_G)$ satisfying  
\begin{eqnarray*}
	&&	0\in \alpha\nabla f(\bar x,\bar y)-\alpha W(\bar x,\bar y, u,v) \times \{0\}+\nabla g(\bar x,\bar y)^T \nu_g
		+\nabla G(\bar x,\bar y)^T \nu_G,\\
		&&	\alpha\geq0,\ 0\leq\nu_g\perp  g(\bar x,\bar y), \nu_g \perp\nabla g(\bar x,\bar y)(u,v),\ 0\leq\nu_G\perp G(\bar x,\bar y), \nu_G \perp \nabla G(\bar x,\bar y)(u,v),
	\end{eqnarray*}
	where $W(\bar x,\bar y, u,v ):=\{\nabla_x f(\bar x,\bar y)+\nabla_x g(\bar x,\bar y)^T \lambda| \lambda \in \Lambda (\bar x,\bar y)\cap \{\nabla g(\bar x,\bar y)(u,v)\}^\perp \}$
	and there exist sequences $t_k\downarrow0,\ (u^k,v^k)\rightarrow(u,v)$ such that (\ref{sequencial1})-(\ref{sequencial3}) hold.
Then there exists a vector $(\lambda^V,\lambda_g,\lambda_G,\lambda )\in \mathbb{R}^{1+p+q+p}$ satisfying 
\begin{align*}
0&=\nabla_x F(\bar x,\bar y)-\lambda^V\nabla_x g(\bar x,\bar y)^T \lambda +\nabla_x g(\bar x,\bar y)^T\lambda_g+\nabla_x G(\bar x,\bar y)^T\lambda_G,\\
0&=\nabla_y F(\bar x,\bar y)+\lambda^V\nabla_y f(\bar x,\bar y)+\nabla_y g(\bar x,\bar y)^T\lambda_g+\nabla_y G(\bar x,\bar y)^T\lambda_G\\
&\lambda^V\geq0,\ 0\leq\lambda_g\perp g(\bar x,\bar y),\ \lambda_g\perp\nabla g(\bar x,\bar y)(u,v)\ ,
	0\leq\lambda_G\perp G(\bar x,\bar y),\ \lambda_G\perp\nabla G(\bar x,\bar y)(u,v)\\
	&0= \nabla_y f(\bar x,\bar y)+\nabla_y g(\bar x,\bar y)^T\lambda, 0\leq\lambda \perp g(\bar x,\bar y), \lambda\perp\nabla g(\bar x,\bar y)(u,v).
\end{align*}
\end{thm}
\beginproof 
{By Proposition \ref{Semic},  since $\mathcal F(x)$ is RCR-regular at $(\bar x,\bar y)$ and $S(x)$ is inner semi-continuous at $(\bar x,\bar y)$ in direction $u$,  $V(x)$ is directional differentiable at $\bar x$ in direction $u$ and $V'(\bar x;u)=\min_{v\in \mathbb{L}(\bar x, \bar y;u)}\nabla f(\bar x, \bar y)(u,v).$ Since $(u,v)\in C(\bar x,\bar y)$, we have $\nabla f(\bar x,\bar y)(u,v)-V'(\bar x;u)=0$.  Hence $v\in\Sigma(\bar x,\bar y,u)= \{v\in \mathbb L( x, y; u)| V'( x;u)=\nabla f(x, y)(u,v)\}$.  Moreover by Theorem \ref{estimates}(ii), $V(x)$ is directionally Lipschitz continuous at $\bar x$ in direction $u$ and by Theorem \ref{subdiff}
$\partial^c V(\bar x; u) \subseteq W(\bar x,\bar y, u,v)$.
The rest of result follows  from Theorem \ref{opt}.}
\endproof

The following example  verifies Theorem \ref{Thm4.4}. 
For this example, $S(x)$ is not inner semi-continuous at $\bar x$ but it is directional inner semi-continuous,  the classical quasi-normality fails but the directional quasi-normality holds.
\begin{example} Consider the following bilevel program
	\begin{eqnarray*}
	\quad & \displaystyle \min_{x,y}& F(x,y):=(\sqrt{3}x-y-\sqrt{3})^2+x+\sqrt{3}y+3\\
		& {\rm s.t.}&  y \in S(x):=\arg\min_y \{1-(x-y)^2: (x-1)^2+y^2-4\leq0, -\sqrt{3}x-y-\sqrt{3}\leq0\}.
	\end{eqnarray*}
	It is easy to verify that 
	\begin{eqnarray}\label{solution}
	S(x)&=&\left\{\begin{array}{ll}
	\sqrt{4-(x-1)^2},\ &-1\leq x<0,\\
	\{-\sqrt{3},\sqrt{3}\},\ &x=0,\\
	-\sqrt{4-(x-1)^2},\ &0<x\leq3.
	\end{array}
	\right. \\
	 V(x)&=&\left\{\begin{array}{ll}
	1-(x-\sqrt{4-(x-1)^2})^2,\ &-1\leq x<0,\\
{-2},\ &x=0,\\
	1-(x+\sqrt{4-(x-1)^2})^2,\ &0<x\leq3.
	\end{array}
	\right.\label{solution2}
	\end{eqnarray} 
	Note that the value function is directionally Lipschitz continuous at $\bar x=0$ but not smooth.
	 The global optimal solution of the bilevel program  is $(\bar x,\bar y)=(0,-\sqrt{3})$.
	By $(\ref{solution})$, $S(x)$ is inner semi-continuous at $\bar y$ in any direction $u>0$. Indeed, for any sequence $x\rightarrow\bar x$ in direction $u>0$, $S(x)\rightarrow \bar y$. It follows that $S(\bar x; u)=\{\bar y\}$.  Note that since for any sequence $x\rightarrow\bar x$ in direction $u<0$, $S(x)\not \rightarrow \bar y$, $S(x)$ is not inner semi-continuous at $\bar x$.

	Denote by $f(x,y):=1-(x-y)^2, g_1(x,y):=(x-1)^2+y^2-4, g_2(x,y):=-\sqrt{3}x-y-\sqrt{3}.$
	Then 
	\begin{equation*}
	\nabla F(\bar x,\bar y)=\left [\begin{matrix}
	1\\
	\sqrt{3}
	\end{matrix}\right ], \quad \nabla f(\bar x,\bar y)=\left [\begin{matrix}
	-2\sqrt{3}\\
	2\sqrt{3}
	\end{matrix}\right ] \quad \nabla g_1(\bar x,\bar y)=\left [\begin{matrix}
	-2\\
	-2\sqrt{3}
	\end{matrix}\right ],  \quad \nabla g_2(\bar x,\bar y)=\left [\begin{matrix}
	-\sqrt{3}\\
	-1
	\end{matrix}\right ].
	\end{equation*}
	It is easy to see that the rank of the gradient vectors $\{\nabla_yg_1(x, y),\nabla_yg_2( x, y)\}$ is always equal to $1$ around $(\bar x,\bar y)$ and hence, RCR-regularity holds at $(\bar x,\bar y)$.
Since $g_1(\bar x,\bar y)=0, g_2(\bar x,\bar y)=0$, 
$$\Lambda(\bar x,\bar y):=\{(\lambda^1, \lambda^2)\in\mathbb R^2_+| 2\sqrt{3}-2\sqrt{3}\lambda^1-\lambda^2=0\}\notag .$$
	Then by Theorem \ref{Thm4.4}, $V(x)$ is directionally Lipschitz continuous and directionally differentiable in direction $u>0$ and 
	\begin{align*}
	V'(\bar x;u)=&\max\{\nabla_x\mathcal L(\bar x,\bar y;\lambda^1,\lambda^2) u: (\lambda^1,\lambda^2)\in\Lambda(\bar x,\bar y)\}\notag\\
	=&\max\{(-2\sqrt{3}-2\lambda^1-\sqrt{3}\lambda^2) u| (\lambda^1,\lambda^2)\in\mathbb R^2_+,2\sqrt{3}-2\sqrt{3}\lambda^1-\lambda^2=0\}\notag\\
	=&\max\{(-2\sqrt{3}+4\lambda^1-6) u| 0\leq\lambda^1\leq1\}\notag\\
	=&-(2\sqrt{3}+2)u.
	\end{align*}
Moreover we can verify that this statement is correct by the expression  (\ref{solution2}).
	Now we prove that the directional quasi-normality holds at $(\bar x,\bar y)$.
The critical cone can be calculated as
	\begin{align*}
	 C(\bar x,\bar y)&:=\{(u,v)|\nabla F(\bar x,\bar y)(u,v)\leq0, \nabla f(\bar x,\bar y)(u,v)-V'(\bar x;u)=0,\nabla g(\bar x,\bar y)(u,v)\leq 0\} \\
	&= \{(u,v)| u+\sqrt{3} v=0, \sqrt{3} u+v\geq 0\}.
	\end{align*}
Let $\bar u=\sqrt{3}$ and $\bar v=-1$, we have  $(\bar u,\bar v)\in C(\bar x,\bar y)$. 
Since $g_1(\bar x,\bar y)=g_2(\bar x,\bar y)=0, \nabla g_1(\bar x,\bar y)(\bar u,\bar v)=0, \nabla g_2(\bar x,\bar y)(\bar u,\bar v)=-\sqrt{3} \bar u-\bar v\not =0$, we have
\begin{eqnarray*}
W(\bar x,\bar y,\bar u,\bar v)&:=&  \{\nabla_x f(\bar x,\bar y)+\nabla_x g(\bar x,\bar y)^T\lambda_g| \lambda_g \in \Lambda(\bar x,\bar y) \cap \{\nabla g(\bar x,\bar y)(\bar u,\bar v)\}^\perp\}\\
&=& \{ -2\sqrt{3}-2 \lambda_g^1| \lambda_g^1\geq 0, 2\sqrt{3}-2 \sqrt{3}\lambda_g^1=0\}\\
&=& \{ -2 \sqrt{3}-2\}.
\end{eqnarray*}
Since $(\bar u,\bar v)\in C(\bar x,\bar y)$, by (\ref{Sigma})  we have $\bar v\in \Sigma (\bar x,\bar y, \bar u)$.  Therefore by Theorem \ref{subdiff}, we have $\partial^c V(\bar x; \bar u)\subseteq W(\bar x,\bar y,\bar u,\bar v)$. Since $V(x)$ is a function of one variable, we can verify by the expression of the value function (\ref{solution2}) that
$$\partial^c V(\bar x; \bar u)=W(\bar x,\bar y,\bar u,\bar v)=\{-2 \sqrt{3}-2\}.$$
Let 
	$\alpha, \nu_1, \nu_2
	$ be such that
	\begin{eqnarray}\label{3x}
	&& 0\in\alpha(\nabla_x f(\bar x,\bar y)-W(\bar x,\bar y, \bar u, \bar v ))+\nu_1\nabla_xg_1(\bar x,\bar y)+\nu_2\nabla_xg_2(\bar x,\bar y), \\
	&&	0=\alpha\nabla_y f(\bar x,\bar y)+\nu_1\nabla_yg_1(\bar x,\bar y)+\nu_2\nabla_yg_2(\bar x,\bar y),\label{3y}\\
	&& {\nu_2\nabla g_2(\bar x,\bar y)(\bar u,\bar v)=0}, \alpha \geq0,\ \nu_1\geq0,\ \nu_2\geq0 \label{3z}
	\end{eqnarray}
	and there exist sequences $t_k\downarrow0,\ (u^k,v^k)\rightarrow (\bar u,\bar v)$, such that
	\begin{eqnarray}
	&&f(\bar x+t_ku^k,\bar y+t_kv^k)-V(\bar x+t_ku^k)>0\ \mbox{if}\ \alpha>0,\label{seq1}\\
	&& g_1(\bar x+t_ku^k,\bar y+t_kv^k)>0\ \mbox{if}\ \nu_1>0.\label{seq4}\\
	&& g_2(\bar x+t_ku^k,\bar y+t_kv^k)>0\ \mbox{if}\ \nu_2>0.\label{seq5}
	\end{eqnarray}
	(\ref{3z}) implies that $\nu_2=0$ and (\ref{seq5}) will not {be} needed.
We now show the conditions (\ref{3x})-(\ref{seq4}) can only hold if $\alpha=\nu_1=\nu_2=0$. By $(\ref{3y})$, $2\sqrt{3}\alpha-2\sqrt{3}\nu_1=0$.  Hence $\alpha=\nu_1$. To the contrary, assume $\alpha>0$. Then $\nu_1=\alpha>0$.
	 Let  $t_k\downarrow0,\ (u^k,v^k)\rightarrow (\bar u,\bar v)$ be arbitrary and suppose that (\ref{seq4}) holds. Then $g_1(x^k,y^k)>0$ for $(x^k,y^k):=(\bar x+t_ku^k,\bar y+t_kv^k)$. It follows that   $y^k<-\sqrt{4-(x^k-1)^2}$. Since  $\nabla_yf(x^k,-\sqrt{4-(x^k-1)^2)}=2(x^k+\sqrt{4-(x^k-1)^2)}>0$ and  $y^k<-\sqrt{4-(x^k-1)^2}$ we have $f(x^k,y^k)<f(x^k,-\sqrt{4-(x^k-1)^2})=V(x^k)$, where the last equality follows from  $(\ref{solution})$. Hence (\ref{seq1}) does not hold.
	The contradiction show that  $(\alpha,\nu_1,\nu_2)=(0,0,0)$ and directional quasi-normality holds at $(\bar x,\bar y)$ in direction $(\bar u,\bar v)$.

	By now, the conditions in Theorem \ref{Thm4.4} are all verified and so the directional KKT condition should hold at $(\bar x,\bar y)$. 
	That is, there exists a nonzero vector $(\lambda_V,\lambda,\lambda_g)\in\mathbb R^{1+2+2}$
	such that
	\begin{eqnarray*}
&&	0=1-{\lambda_V}(-2\lambda^1-\sqrt{3}\lambda^2)-2\lambda^1_g-\sqrt{3}\lambda^2_g,\\
&&	0=\sqrt{3}+{\lambda_V}2\sqrt{3}-2\sqrt{3}\lambda^1_g-\lambda^2_g,\\
	&& \lambda_g, \lambda \in \Lambda(\bar x,\bar y), \lambda_g \perp \nabla g(\bar x,\bar y)(\bar u,\bar v), \lambda \perp \nabla g(\bar x,\bar y)(\bar u,\bar v).
	\end{eqnarray*}
	Obviously the vectors
	$(\lambda_V,\lambda,\lambda_g):=(\frac{1}{2},(1,0),(1,0))$ satisfies the above conditions.
	
{As we have mentioned before, NNAMCQ and FOSCMS  always fail for (BP). In this example, the quasi-normality also fails at $(\bar x,\bar y)$. Indeed, let $(\alpha,\nu_1,\nu_2)=(1,1,0)$. We have $(\alpha,\nu_1,\nu_2)$ satisfies $(\ref{3x})$ and $(\ref{3y})$. And choose $(x^k,y^k):=(-1/k-\sqrt{4-(1/k+1)^2}-1/k)$, which converges to $(\bar x,\bar y)$. By $(\ref{solution})$, we have 
\begin{align*}
&f(x^k,y^k)=1-\left(\sqrt{4-\left(1/k+1\right)^2}\right)^2>1-\left(1/k+\sqrt{4-\left(1/k+1\right)^2}\right)^2=V(x^k),\\
&g_1(x^k,y^k)=\left(1/k+1\right)^2+\left(\sqrt{4-\left(1/k+1\right)^2}+1/k\right)^2-4>0,
\end{align*}
By the definition of the classical quasi-normality defined in \cite[Definition 4.2]{GYZ} (one can refer to Definition \ref{qp} for the case $u=0$), this means that the quasi-normality fails at $(\bar x,\bar y)$.}
\end{example}

\section*{Acknowlegement} The authors would like to thank the anonymous referees for their helpful suggestions and comments.


\end{document}